\documentclass[10pt, reqno]{amsart}
\usepackage{amsmath,amssymb,amscd,mathrsfs,amscd}
\usepackage{graphics,verbatim}
\newdimen\hoogte    \hoogte=12pt    
\newdimen\breedte   \breedte=14pt   
\newdimen\dikte     \dikte=0.5pt    

\newenvironment{Young}{\begingroup
       \def\vr{\vrule height0.8\hoogte width\dikte depth 0.2\hoogte}
       \def\fbox##1{\vbox{\offinterlineskip
                    \hrule height\dikte
                    \hbox to \breedte{\vr\hfill##1\hfill\vr}
                    \hrule height\dikte}}
       \vbox\bgroup \offinterlineskip \tabskip=-\dikte \lineskip=-\dikte
            \halign\bgroup &\fbox{##\unskip}\unskip  \crcr }
       {\egroup\egroup\endgroup}
\def\diagram#1{\relax\ifmmode\vcenter{\,\begin{Young}#1\end{Young}\,}\else%
              $\vcenter{\,\begin{Young}#1\end{Young}\,}$\fi}

\theoremstyle{plain}
\newtheorem{thm}{Theorem}[subsection]
\newtheorem{lem}[thm]{Lemma}

\newtheorem{prp}[thm]{Proposition}

\newtheorem{cor}[thm]{Corollary}

\newtheorem{rmk}[thm]{Remark}

\newtheorem{exa}[thm]{Example}

\newenvironment{pff}{{\em Proof:}}{\QED}

\numberwithin{equation}{subsection}

\newcommand{\af}{\alpha}
\newcommand{\bt}{\beta}
\newcommand{\gm}{\gamma}
\newcommand{\dt}{\delta}

\newcommand{\io}{\iota}

\newcommand{\sm}{\sigma}
\newcommand{\kp}{\kappa}

\newcommand{\Gm}{\Gamma}
\newcommand{\Dt}{\Delta}

\newcommand{\del}{\partial}

\newcommand{\Q}{{\mathbb{Q}}}
\newcommand{\Z}{{\mathbb{Z}}}

\newcommand{\C}{{\mathbb{C}}}

\newcommand{\bbF}{{\mathbb{F}}}

\newcommand{\Sl}{{\mathfrak{sl}}}

\newcommand{\g}{{\mathfrak{g}}}
\newcommand{\h}{{\mathfrak{h}}}
\newcommand{\n}{{\mathfrak{n}}}

\newcommand{\q}{{\mathfrak{q}}}

\renewcommand{\b}{{\mathfrak{b}}}

\newcommand{\A}{{\mathcal{A}}}
\newcommand{\HC}{{\mathcal{HC}}}
\renewcommand{\P}{{\mathcal{P}}}
\renewcommand{\L}{{\mathcal{L}}}
\newcommand{\M}{{\mathcal{M}}}

\newcommand{\U}{{\mathcal{U}}}
\newcommand{\curlyQ}{{\mathcal{Q}}}
\newcommand{\F}{{\mathcal{F}}}
\newcommand{\W}{{\mathcal{W}}}

\newcommand{\G}{{\mathcal{G}}}
\newcommand{\GL}{{\mathcal{GL}}}

\newcommand{\ch}{{\operatorname{ch}\,}}
\newcommand{\ind}{{\operatorname{Ind}}}
\newcommand{\res}{{\operatorname{Res}}}

\newcommand{\Hom}{{\operatorname{Hom}}}

\newcommand{\End}{{\operatorname{End}}}

\newcommand{\Rep}{{\operatorname{Rep}}}
\newcommand{\Proj}{{\operatorname{Proj}}}

\newcommand{\height}{{\operatorname{ht}}}

\newcommand{\cont}{{\operatorname{cont}}}
\newcommand{\gdim}{{\dim_q\,}}

\newcommand{\ui}{{\underline{i}}}
\newcommand{\uj}{{\underline{j}}}

\newcommand{\la}{{\langle}}
\newcommand{\ra}{{\rangle}}

\newcommand{\va}{{\mathbf{1}}}

\newcommand\dubc{, \ldots ,}

\newcommand{\andeqn}{\,\,\,\,\,\, {\mbox{and}} \,\,\,\,\,\,}
\newcommand{\QED}{\rule{0.4em}{2ex}}

\input xy
\xyoption{all}

\begin{document}
\title{Representations of Quiver Hecke Algebras via Lyndon Bases}

\author{David Hill}
\address{Department of Mathematics \\
            University of California, Berkeley \\
            Berkeley, CA 94720-3840}
\email{dhill1@math.berkeley.edu}

\author{George Melvin}
\address{Department of Mathematics \\
            University of California, Berkeley \\
            Berkeley, CA 94720-3840}
\email{gmelvin@math.berkeley.edu}

\author{Damien Mondragon}
\address{Department of Mathematics \\
            University of California, Berkeley \\
            Berkeley, CA 94720-3840}
\email{damien@math.berkeley.edu}

\thanks{Research of the first author was partially supported by
NSF EMSW21-RTG grant DMS-0354321}\

\begin{abstract} A new class of algebras has been introduced by Khovanov and Lauda and independently by Rouquier.
These algebras categorify one-half of the Quantum group associated to arbitrary
Cartan data. In this paper, we use the combinatorics of Lyndon words to
construct the irreducible representations of those algebras associated to
Cartan data of finite type. This completes the classification of simple modules for the quiver Hecke algebra initiated by Kleshchev and Ram.
\end{abstract}

\maketitle

\section{Introduction}

\subsection{} Recently, Khovanov and Lauda \cite{khl1,khl2} and Rouquier
\cite{rq} have independently introduced a remarkable family of graded algebras,
$H(\Gm)$, defined in terms of quivers associated to the Dynkin diagram, $\Gm$, of a
symmetrizable Kac-Moody algebra, $\g$. These algebras have been given
several names, including Khovanov-Lauda-Rouquier algebras, quiver
nil-Hecke algebras, quiver Hecke algebras, and ``the rings $R(\nu)$'' (here $\nu$ refers to an element in the positive cone $\curlyQ^+$ inside the root lattice of $\g$). The
main property of these algebras is that
$$K(\Gm)\cong \U_\A^*(\n)$$ as \emph{twisted} bimodules, where $K(\Gm)$ is the Grothendieck group of
the full subcategory, $\Rep(\Gm)$, of finite dimensional graded
$H(\Gm)$-modules, $\n$ is a maximal nilpotent subalgebra of $\g$, and
$\U_\A^*(\n)$ is an integral form of the quantized enveloping algebra,
$\U_q(\n)$.

Further evidence of the importance of these algebras was obtained
in \cite{bk}. In this work, Brundan and Kleshchev showed that when $\Gm$ is of
type $A_\infty$ or $A_{\ell-1}^{(1)}$, there is an isomorphism between blocks
of cyclotomic Hecke algebras of symmetric groups, and blocks of a
corresponding quotient of $H(\Gm)$. Moreover, this isomorphism applies equally
well to the Hecke algebra and its rational degeneration, depending only on
$\Gm$ and the underlying ground field.
In light of the work \cite{bk2}, it is expected that a similar relationship
should hold between interesting quotients of $H(\Gm)$ and cyclotomic Hecke-Clifford algebras when $\Gm$ is of type $B_\infty$ and $A_{2\ell}^{(2)}$. For these reasons, we choose to use the name ``\emph{quiver Hecke algebra}'' to describe $H(\Gm)$.

\subsection{} In \cite{hks}, Hill, Kujawa and Sussan investigated the
representation theory of the (degenerate) affine Hecke-Clifford algebra, $\HC(d)$,
over $\C$. In this paper, the authors constructed an analogue of the
Arakawa-Suzuki functor \cite{as} between the category $\mathcal{O}$ for the Lie
superalgebra $\q(n)$ and a certain category, $\Rep\HC(d)$, of \emph{integral} finite dimensional
modules for $\HC(d)$. By considering small rank instances of the functor, the
authors obtained analogues of Zelevinsky's \emph{segment representations},
\cite{bz,z}, for $\HC(d)$. More generally, the Verma modules for $\q(n)$
correspond under the functor to certain induced modules, which by \cite[Theorem
4.4.10]{hks} have unique irreducible quotients. The authors went on to obtain a
construction of all the irreducible integral representations using the
combinatorics of Lyndon words together with \cite[Theorem 7.17]{bk2}.

It is instructive to describe their result in more detail. To this end, let $\Gm$ be a Dynkin diagram of finite type with nodes labelled by the index set $I$, fix a total ordering, $\leq$, on $I$. Let $\F$ be the free associative algebra generated by the letters $[i]$, $i\in I$, with the  concatenation product $[i_1]\cdots[i_k]=[i_1,\ldots,i_k]$, and give the monomials in $\F$ the lexicographic ordering determined by $I$. It was first notices in \cite{lr} that certain monomials in $\F$ associated to this ordering, called \emph{good Lyndon words}, and their non-increasing products, called \emph{good words}, naturally determine various bases of the quantized enveloping algebra, $\U_q(\n)$, of a maximal nilpotent subalgebra of the semisimple Lie algebra $\g$ associated to $\Gm$. This observation was further developed in the prophetic paper of Leclerc, \cite{lec}, where it was first suggested that the bases arising from these combinatorics should naturally correspond to representations of affine Hecke algebras, cf. \cite[Sections 6-7]{lec}, specifically \cite[Theorem 47, Conjecture 52]{lec}.

In \cite{hks}, the authors noticed that the character of each segment representation of $\HC(d)$ corresponds in a natural way to a dual canonical basis element labeled by a good Lyndon word in type $B$ with respect to the standard Dynkin ordering on $I$ (specialized at $q=1$), see \cite[Proposition 4.1.3, Theorem 4.1.8, Proposition 8.2.12]{hks}. This was a nontrivial observation since it applied only after redeveloping the theory so that monomials are ordered lexicographically from \emph{right-to-left}, a technicality imposed by the functor, cf. \cite[Lemma 8.2.13]{hks}. This choice had the effect of drastically simplifying both the good Lyndon words in type $B$, and their associated dual canonical basis elements. More generally, the characters of standard modules naturally correspond to dual PBW basis elements labeled by good words (again, at $q=1$), \cite[Theorem 8.5.1]{hks}. Finally, applying \cite[Theorem 4.4.10]{hks} completed the construction, \cite[Theorem 8.5.5]{hks}.

Motivated by the results of \cite{hks} and the conjectured connection between $\HC(d)$ and quiver Hecke algebras of type $B$, we initiated a study of the representation theory of the category $\Rep(\Gm)$, for $\Gm$ of classical finite type, using the combinatorics of Lyndon words with respect to the standard Dynkin ordering on $I$ and the right-to-left lexicographic ordering described in \cite{hks}. Indeed, we first observed that this simplified the good Lyndon words in every type (except for the long roots in type $C$, which remain the same). Subsequently, we worked out the corresponding dual canonical basis elements, $b^*_l$, associated to each good Lyndon word, $l$, and constructed representations, $\va_l$, with character $b^*_l$. The standard module, $\M(g)$, associated to a good word $g$ is the module obtained by parabolic induction: $$\M(g)=\ind\, \va_{l_1}\boxtimes\cdots\boxtimes\va_{l_k}\{c_g\},$$ where $g=l_1\cdots l_k$ is the canonical factorization of $g$ as a non-increasing product of good Lyndon words and the term $\{c_g\}$ refers to a grading shift. These standard modules have the property that their characters are given by dual PBW basis elements labelled by the corresponding good word, and, therefore,  give a basis for the Grothendieck group, $K(\Gm)$.

\subsection{} While this paper was in production, Kleshchev and Ram completed their own investigation of $\Rep(\Gm)$ using the combinatorics of Lyndon words, for $\Gm$ of \emph{arbitrary finite type}. To describe this paper in more detail, give $I$ an arbitrary total ordering. The authors called an irreducible $H(\Gm)$-module \emph{cuspidal} if its character is given by a dual canonical basis element associated to a good Lyndon word, cf. \cite[Lemma 6.4]{klram2}. They went on to prove an amazing lemma. Namely, given a cuspidal representation, $\va_l$, the module $$\M(l^k)=\ind\, \underbrace{\va_{l}\boxtimes\cdots\boxtimes\va_l}_{k\mbox{ times}}\,\{c_{l^k}\}$$  remains irreducible for all $k>1$, \cite[Lemma 6.6]{klram2}. We want to point out that this lemma applies equally well to all possible orderings on $I$ and all finite root systems. Combining \cite[Lemma 6.6]{klram2} with a straightforward Frobenius reciprocity argument shows that the standard module $\M(g)$ has a unique irreducible quotient $\L(g)$, \cite[Theorem 7.2]{klram2}. In this way, Kleshchev and Ram reduced the study of $\Rep(\Gm)$ to the construction of cuspidal representations. They went on to construct all cuspidal representations in types $ABCDG$ as well as $E_6$ and $E_7$ using the standard Dynkin ordering, the good Lyndon words in \cite{lr,lec}, and the corresponding root vectors in \cite[Section 8]{lec}, cf. \cite[Section 8]{klram2}. In type $A$ they produced cuspidal representations for all orderings on $I$.

\subsection{} Given the beautiful results in \cite{klram2}, we expanded the goal of this paper. In particular, our main result is a complete determination of the cuspidal representations of $H(\Gm)$ in all finite types using our ordering, Theorem \ref{T:CuspidalsExist}. We would like to point out several advantages of our approach. First, in classical type, our cuspidal representations tend to be much simpler than those appearing in \cite{klram2}. More specifically, in types $BCD$, our representations generally have dimension at most 2 (with the exception of the long roots in type $C$). In contrast, the cuspidal modules constructed by Kleshchev an Ram generally have dimensions that grow with the height of an associated positive root. Another advantage can be seen when considering the case of $E_8$. The main difficulty for Kleshchev and Ram is that not all the $E_8$ root vectors are homogeneous in the sense of \cite{klram}, see \cite[Section 8]{klram2}. On the other hand, in our ordering, all good Lyndon words in type $E_8$ are homogeneous. Finally, in $\S$\ref{SS:antiautomorphism of W}, \ref{SS:Automorphism}, and \ref{SS:Twisting} we explain exactly how to relate the right-to-left lexicographic ordering used here to the more standard left-to-right lexicographic ordering in \cite{lec} and \cite{klram2}.

In this paper, we only use the half of the bialgebra structure of $K(\Gm)$ coming from parabolic induction. It would also be interesting to consider the structure coming from restriction and compare the work here to that of Lauda and Vazirani, \cite{vl}.

Finally, we would like to point out that the description of the simple modules for the quiver Hecke algebra of type $B$ is nearly identical to the description of the irreducible $\HC(d)$-modules appearing in \cite{hks}. In particular, it is possible to define an action of (an appropriately defined) \emph{quiver Hecke-Clifford} superalgebra of type $B$ on the segment representations of $\HC(d)$. Moreover, this action extends easily to standard modules. Based on small rank calculations, we conjecture that this action factors through the unique simple quotients. We feel that an investigation of this phenomenon should shed light into the relationship between the type $B$ quiver Hecke algebra and the Hecke-Clifford algebra, but this is a topic of another paper.

\subsection{} The remainder of the paper is organized as follows. In Section \ref{S:QuantGrps} we describe the embedding of the quantum group $\U_q(\n)$ inside the $q$-shuffle algebra $\F$ and describe the combinatorics of Lyndon words in our set-up following \cite{lec} and \cite[Section 8]{hks} closely. In Section \ref{S:QHA} we introduce the quiver Hecke algebra and describe some of the basic properties of the category $\Rep(\Gm)$. In Section \ref{S:StdRepsandSimples} we introduce cuspidal representations, and standard representations and state the main theorem of the paper, Theorem \ref{T:CuspidalsExist}. In Section \ref{S:IdentGoodLyndon} we determine the good Lyndon words and corresponding root vectors, and Section \ref{A:SegReps} contains the construction of cuspidal representations. Finally,  Appendix \ref{S:Appendix} contains the calculations relevant to Section \ref{S:IdentGoodLyndon}.

\subsection*{Acknowledgement} We would like to thank both Alexander Kleshchev and Arun Ram for encouraging us to work out the cuspidal representations in types $E$ and $F$, as well as for their extremely useful comments on an earlier draft of the paper. The first author would additionally like to thank the algebra group in the department of Mathematics at the University of California, Berkeley, and particularly his sponsor, Mark Haiman, for giving him the opportunity to teach a graduate course in the spring of 2009, where the idea to write this paper was first realized.

\section{Quantum Groups}\label{S:QuantGrps}

\subsection{Root Data}\label{SS:Root Data}
Let $\g$ be a simple finite dimensional Lie algebra of rank $r$ over $\C$, with
Dynkin diagram $\Gm$ and let $I$ denote the set of labels of the nodes of
$\Gm$. Let $\U_q(\g)$ be the corresponding quantum group over $\Q(q)$ with
Chevalley geneators $e_i,f_i$, $i\in I$. Let $\n\subseteq\g$ be the
subalgebra generated by the $e_i$, $i\in I$. Let $\Dt$ be the root system of
$\g$ relative to this decomposition, $\Dt^+$ the positive roots, and
$\Pi=\{\af_i|i\in I\}$ the simple roots. Let $\curlyQ$ be the root lattice and
$\curlyQ^+=\sum_{i\in I}\Z_{\geq 0}\af_i$. Let $A=(a_{ij})_{i,j\in
I}$ be the Cartan matrix of $\g$  and $(\cdot,\cdot)$ denote
symmetric bilinear form on $\h^*$ satisfying
\[
a_{ij}=\frac{2(\af_i,\af_j)}{(\af_i,\af_i)},\;\;\;
d_i=\frac{(\af_i,\af_i)}2\in\{1,2,3\}.
\]
Let $q_i=q^{d_i}$. Define the $q$-integers and $q$-binomial
coefficients:
\[
[k]_i=\frac{q_i^k-q_i^{-k}}{q_i-q_i^{-1}},\;\;\;[k]_i!=[ k]_i\cdots[ 2]_i[ 1]_i,\;\;\; \left[{m \atop k}\right]_i
=\frac{[m]_i!}{[ k]_i![m-k]_i!}.
\]

For later purposes, we also define the following. Let $\nu\in\curlyQ^+$, say $\nu=\sum_{i\in I}c_i\af_i$. Define the
\emph{height} of $\nu$:
\[
\height(\nu)=\sum_{i\in I}c_i.
\]
Next, given $\ui=(i_1,\ldots,i_d)\in I^d$, define the \emph{content} of $\ui$
by
\[
\cont(\ui)=\sum_{i\in I}n_i\af_i,\;\;\;n_i=\#\{j=1,\ldots,d\,|\,i_j=i\}.
\]
Finally, if $\height(\nu)=d$, set $I^\nu=\{\ui\in I^d|\cont(\ui)=\nu\}$. Let
$S_d$ denote the symmetric group on $d$ letters, generated by simple
transpositions $s_1,\ldots,s_{d-1}$. Then, $S_d$ acts $I^d$ by place
permutation and we denote this action by $w\cdot\ui$, $w\in S_d$, $\ui\in I^d$.
Observe that the orbits of this action are precisely the sets $I^\nu$ with
$\height(\nu)=d$.

\subsection{Embedding of $\U_q(\n)$ in the Quantum Shuffle
Algebra}\label{SS:Embedding}

The algebra $\U_q:=\U_q(\n)$ is a quotient of the free algebra generated by
the Chevalley generators $e_i$, $i\in I$ by the relations
\[
\sum_{r+s=1-a_{ij}}(-1)^r\left[{1-a_{ij}\atop r}\right]_ie_i^re_je_i^s=0.
\]
It is naturally $\curlyQ^+$-graded by assigning to $e_i$ the degree $\af_i$.
Let $|u|$ be the $\curlyQ^+$-degree of a homogeneous element $u\in\U_q$.

In \cite{ka}, Kashiwara proved that there exist $q$-derivations $e_i'$, $i\in
I$ given by
\[
e_i'(e_j)=\dt_{ij}\andeqn e_i'(uv)=e_i'(u)v+q^{-(\af_i,|u|)}ue_i'(v)
\]
for all homogeneous $u,v\in\U_q$. For each $i\in I$, $e_i'(u)=0$ if, and only
if $|u|=0$.

Now, let $\F$ be the free associative algebra over $\Q(q)$ generated by the set
of letters $\{[i]|i\in I\}$. Letters should not be confused with q-integers, which always
occur with a subscript. Write
$[i_1,\ldots,i_k]:=[i_1]\cdot[i_2]\cdots[i_k]$, and let $[]$ denote the empty
word. The algebra $\F$ is $\curlyQ^+$ graded by assigning the degree $\af_i$ to
$[i]$ (as before, let $|f|$ denote the $\curlyQ^+$-degree of a homogeneous
$f\in\F$). Notice that $\F$ also has a \emph{principal grading} obtained by
setting the degree of a letter $[i]$ to be 1; let $\F_d$ be the $d$th graded
component in this grading.

Now, define the (quantum) shuffle product, $*$, on $\F$ inductively by
\begin{align}\label{E:inductiveqshuffle}
(x\cdot[i])*(y\cdot[j])=(x*(y\cdot[j])\cdot[i]+q^{-(|x|+\af_i,\af_j)}((x\cdot[i])*y)\cdot[j],\;\;\;x*[]=[]*x=x.
\end{align}
Iterating this formula yields
\begin{align}\label{E:Shuffle}
[i_1,\ldots,i_\ell]*[i_{\ell+1},\ldots,i_{\ell+k}] =\sum_{w\in
D_{(\ell,k)}}q^{-e(w)}[i_{w^{-1}(1)},\ldots,i_{w^{-1}(k+\ell)}]
\end{align}
where $D_{(\ell,k)}$ is the set of minimal coset representatives in $S_{\ell+k}/S_\ell\times S_k$ and
\[
e(w)=\sum_{\substack{s\leq\ell<t\\w(s)<w(t)}}(\af_{i_{s}},\af_{i_{t}}),
\]
see \cite[$\S2.5$]{lec}. The product $*$ is associative and,
\cite[Proposition 1]{lec},
\begin{eqnarray}\label{E:qShuffle}
x*y=q^{-(|x|,|y|)}y\overline{*}x
\end{eqnarray}
where $\overline{*}$ is obtained by replacing $q$ with $q^{-1}$ in the
definition of $*$.

Now, to $f=[i_1,\ldots,i_k]\in\F$, associate $\del_f=e_{i_1}'\cdots
e_{i_k}'\in\End \U_q$, and $\del_{[]}=\operatorname{Id}_{\U_q}$. Then,

\begin{prp}\cite{ro1,ro2,grn} There exists an injective $\Q(q)$-linear homomorphism
\[
\Psi:\U_q\rightarrow(\F,*)
\]
defined on homogeneous $u\in\U_q$ by the formula $\Psi(u)=\sum\del_f(u)f$,
where the sum is over all monomials $f\in\F$ such that $|f|=|u|$.
\end{prp}

Therefore $\U_q$ is isomorphic to the subalgebra $\W\subseteq(\F,*)$
generated by the letters $[i]$, $i\in I$.

Let $\A=\Q[q,q^{-1}]$, and let $\U_\A$ denote the $\A$-subalgebra of $\U_q$
generated by the divided powers $e_i^k/[k]_i!$ ($i\in I$, $k\in\Z_{\geq0}$).
Let $(\cdot,\cdot)_K:\U_q\times\U_q\rightarrow\Q(q)$ denote the unique
symmetric bilinear form satisfying
\begin{align}\label{E:SymmetricForm}
(1,1)_K=1\andeqn(e_i'(u),v)_k=(u,e_iv)_K
\end{align}
for all $i\in I$, and $u,v\in\U_q$. Let
\begin{align}\label{E:DualEnvelope}
\U_\A^*=\{\,u\in\U_q \mid (u,v)_K\in\A\mbox{ for all }v\in\U_\A\,\}
\end{align}
and let $u^*\in\U_\A^*$ denote the dual to $u\in\U_\A$ relative to $(\cdot,\cdot)_K$. It is well known that for $u\in\U_{\A,\nu}$, the map $u^*\mapsto(u^*,?)_K$ defines an isomorphism $\U_{\A,\nu}^*\cong\Hom_\A(\U_{\A,\nu},\A)$.

\begin{rmk}\label{R:LusztigForm} Observe that the form we are using differs slightly from Lustig's bilinear form $(\cdot,\cdot)_L$. They are related by the formula
$$(u,v)_L=\prod_{i\in I}\frac{1}{(1-q_i^2)^{c_i}}(u,v)_K,$$
if $|u|=|v|=\sum_ic_i\af_i$. In particular, if $B$ is a basis of $\U_q$ consisting of homogeneous vectors, then the adjoint basis of $B$ with respect to $(\cdot,\cdot)_K$ and $(\cdot,\cdot)_L$ differ only by some normalization factors. In particular, $B$ is orthogonal with respect to $(\cdot,\cdot)_K$ if, and only if it is orthogonal with respect to $(\cdot,\cdot)_L$.

Throughout this paper, we will shall follow Leclerc and use the form $(\cdot,\cdot)_K$. In $\S$\ref{SS:ModulesandChar} we will explain how both forms arise in representation theory, cf. Example \ref{Exa:KKLform} and Lemma \ref{L:IdentificationofLattices}.
\end{rmk}

Now, given a monomial
\[
[i_1^{a_1},i_2^{a_2},\ldots,i_k^{a_k}]
    =[\underbrace{i_1,\ldots,i_1}_{a_1},\underbrace{i_2,\ldots,i_2}_{a_2},
    \ldots,\underbrace{i_k,\ldots,i_k}_{a_k}]
\]
 with $i_j\neq i_{j+1}$ for $1\leq j<k$, let
 $c_{i_1,\ldots,i_k}^{a_1,\ldots,a_k}=[a_1]_{i_1}!\cdots[a_k]_{i_k}!$, so that
 $(c_{i_1,\ldots,i_k}^{a_1,\ldots,a_k})^{-1}e_{i_1}^{a_1}\cdots e_{i_k}^{a_k}$ is a
 product of divided powers. Let
\[
\F^*_\A=\bigoplus\A c_{i_1,\ldots,i_k}^{a_1,\ldots,a_k}
[i_1^{a_1},i_2^{a_2},\ldots,i_k^{a_k}]
\]
and $\W^*_\A=\W\cap\F^*_\A$. It is known that $\W_\A^*=\Psi(\U_\A^*)$,
\cite[Lemma 8]{lec}.

We close this section by describing some simple involutions of $\F$ which correspond,
on restriction to $\W$, to important involutions on $\U_q$. To this end, for $\nu=\sum_ic_i\af_i\in\curlyQ^+$, define
\begin{align}\label{E:Nfcn}
N(\nu)=\frac12\left((\nu,\nu)-\sum_{i=1}^rc_i(\af_i,\af_i)\right).
\end{align}

\begin{prp}\label{P:BarInv}\cite[Proposition 6]{lec} Let $f=[i_1,\ldots,i_k]$, $|f|=\nu$. Then,

(i) Let $\tau:\F\to\F$
be the $\Q(q)$-linear map defined by $\tau(f)=[i_k,\ldots,i_1]$.
Then, $\tau(x*y)=\tau(y)*\tau(x)$ for all $x,y\in\F$. Hence,
$\tau(\Psi(u))=\Psi(\tau(u))$, where $\tau:\U_q\to\U_q$ is the
$\Q(q)$-linear anti-automorphism which fixes the generators $e_i$.

(ii) Let $-:\F\rightarrow\F$ be the
$\Q$-linear map defined by $\bar{q}=q^{-1}$ and
\[
\overline{f}
=q^{N(\nu)}[i_k,\ldots,i_1].
\]
Then, $\overline{x*y}=\overline{x}*\overline{y}$ for all $x,y\in\F$. Hence,
$\overline{\Psi(u)}=\Psi(\overline{u})$, where $-$ is the bar involution on
$\U_q$.

(iii) Let $\sm:\F\to\F$ be the $\Q$-linear map such that $\sm(q)=q^{-1}$ and
$$\sm(f)=q^{N(\nu)}f.$$ Then, $\sm(x)=\overline{\tau(x)}$
for all $x\in\F$. Hence, $\Psi(\sm(u))=\sm(\Psi(u))$, where
$\sm:\U_q\to\U_q$ is the $\Q$-linear anti-automorphism which sends $q$ to
$q^{-1}$ and fixes the Chevalley generators $e_i$.
\end{prp}

\subsection{Good Words and Lyndon Words}\label{SS:LyndonWords} In what follows,
our conventions differ from those in \cite{lec}. In particular, we order
monomials in $\F$ lexicographically reading from \emph{right to left}. Except
for the type $A$ case, this convention leads to some significant differences in
the good Lyndon words that appear. For the convenience of the reader, we
include $\S$\ref{SS:antiautomorphism of W} which explains the connection between
the combinatorics developed using this ordering to those which arise using the
more common \emph{left to right} lexicographic ordering.

The next two sections parallel \cite[Sections 3,4]{lec} with the statements
of the relevant propositions adjusted to conform to our choice of ordering.

For the remainder of the section, fix an ordering on the set of letters
$\{[i]|i\in I\}$ in $\F$, denoted $\leq$, and order $\Pi$ accordingly. Give the
set of monomials in $\F$ the associated lexicographic order read from right to
left, also denoted $\leq$. That is, set $[i]<[]$ for all $i\in I$ and
\[
[i_1,\ldots,i_k]<[j_1,\ldots,j_\ell]\mbox{ if }i_k<j_\ell,\mbox{ or for some
}m, i_{k-m}<j_{\ell-m}\mbox{ and }i_{k-s}=j_{\ell-s}\mbox{ for all }s<m.
\]
Note that since the empty word is larger than any letter, every word is smaller
than all of its right factors:
\begin{align}\label{E:rightfactors}
[i_1,\ldots,i_k]<[i_j,\ldots,i_k],\mbox{ for all }1<j\leq k.
\end{align}
(For those familiar with the theory, this definition is needed to ensure that
the induced Lyndon ordering on positive roots is convex, cf.
$\S$\ref{SS:PBWandCanonical} below.)

For a homogeneous element $f\in\F$, let $\min(f)$ be the smallest monomial
occurring in the expansion of $f$. A monomial $[i_1,\ldots,i_k]$ is called a
\emph{lower good word} if there exists a homogeneous $w\in\W$ such that
$[i_1,\ldots,i_k]=\min(w)$, and we say that it is \emph{Lyndon on the right} if it is larger
than any of its proper left factors:
\[
[i_1,\ldots,i_j]<[i_1,\ldots,i_k],\mbox{ for any }1\leq j<k.
\]
Except for $\S$\ref{SS:antiautomorphism of W}, we refer to these special words
simply as \emph{good} and \emph{Lyndon}. Let $\G$ denote the set of good words,
$\L$ the set of Lyndon words, and $\GL=\L\cap\G\subset\G$ the set of good
Lyndon words. Also, let $\GL_d\subset\G_d\subset\F_d$ denote the degree $d$
components of $\GL$ and $\G$ in the principal grading. Finally, for $\nu\in\curlyQ^+$, let $\GL_\nu\subset\G_\nu\subset\F_\nu$ be the homogeneous components of $\GL$ and $\G$ in the $\curlyQ^+$ grading.

\begin{lem}\label{L:GoodFactors}\cite[Lemma 13]{lec} Every factor of a good word is
good.
\end{lem}

Because of our ordering conventions, \cite[Lemma 15, Proposition 16]{lec}
become

\begin{lem}\cite[Lemma 15]{lec} Let $l\in\L$, $w$ a monomial such that $w\geq l$. Then, $\min(w*l)=wl$.
\end{lem}
\noindent and
\begin{prp}\label{P:GLproduct}\cite[Proposition 16]{lec} Let $l\in\GL$, and $g\in\G$ with $g\geq l$.
Then $gl\in\G$.
\end{prp}

Hence, we deduce from Lemma \ref{L:GoodFactors} and Proposition
\ref{P:GLproduct} \cite[Proposition 17]{lec}:

\begin{prp}\cite{lr,lec} A monomial $g$ is a good word if, and only if,
there exist good Lyndon words $l_1\geq\ldots\geq l_k$ such that
\[
g=l_1l_2\cdots l_k.
\]
\end{prp}

As in \cite{lec}, we have

\begin{prp}\cite{lr,lec} The map $l\rightarrow|l|$ is a bijection
$\GL\rightarrow\Dt^+$.
\end{prp}

Given $\bt\in\Dt^+$, let $\bt\rightarrow l(\bt)$ be the inverse of the above
bijection (called the Lyndon covering of $\Dt^+$).

We now define the \emph{bracketing} of Lyndon words, that gives rise to the
\emph{Lyndon basis} of $\W$. To this end, given $l\in\L$ such that $l$ is not a
letter, define the standard factorization of $l$ to be $l=l_1l_2$ where
$l_2\in\L$ is a proper left factor of maximal length. Define the $q$-bracket
\begin{align}\label{E:qbracket}
[f_1,f_2]_q=f_1f_2-q^{(|f_1|,|f_2|)}f_2f_1
\end{align}
for homogeneous $f_1,f_2\in\F$ in the $\curlyQ^+$-grading. Then, the bracketing
$\la l\ra$ of $l\in\L$ is defined inductively by $\la l\ra=l$ if $l$ is a
letter, and
\begin{align}\label{E:Lyndonbracketing}
\la l\ra=[\la l_1\ra,\la l_2\ra]_q
\end{align}
if $l=l_1l_2$ is the standard factorization of $l$.

\begin{exa} For $\g$ of type $B_r$ with $I$ given in Table \ref{Tbl:Dynkin} below, we have
\begin{enumerate}
\item $\la [0]\ra=[0]$;
\item $\la [12]\ra=[[1],[2]]_q=[12]-q^{-2}[21]$;
\item $\la[012]\ra=[[0],[12]-q^{-2}[21]]_q=[012]-q^{-2}[021]-q^{-2}[120]+q^{-4}[210]$.
\end{enumerate}
\end{exa}

As is suggested in this example, we have

\begin{prp}\label{P:bracketingtriangularity}\cite[Proposition 19]{lec} For
$l\in\L$, $\la l\ra=l+r$ where $r$ is a linear combination of words $w$ such
that $|w|=|l|$ and $w<l$.
\end{prp}

Any word $w\in\F$ has a canonical factorization $w=l_1\cdots l_k$ such that
$l_1,\ldots,l_k\in\L$ and $l_1\geq\cdots\geq l_k$. We define the bracketing of
an arbitrary word $w$ in terms of this factorization: $\la w\ra=\la
l_1\ra\cdots\la l_k\ra$. Define a homomorphism $\Xi:(\F,\cdot)\to(\F,*)$ by
$\Xi([i])=[i]$. Then,
$\Xi([i_1,\ldots,i_k])=[i_1]*\cdots*[i_k]=\Psi(e_{i_1}\cdots e_{i_k})$. In
particular, $\Xi(\F)=\W$. We have the following characterization of good words:

\begin{lem}\cite[Lemma 21]{lec} The word $w$ is good if and only if it cannot
be expressed modulo $\ker\Xi$ as a linear combination of words $v<w$.
\end{lem}

For $g\in\G$, set $r_g=\Xi(\la g\ra)$. Then, we have

\begin{thm}\label{T:Lyndonbasis}\cite[Propostion 22, Theorem 23]{lec}
Let $g\in\G$ and $g=l_1\cdots l_k$ be the canonical factorization of $g$ as a
nonincreasing product of good Lyndon words. Then
\begin{enumerate}
\item $r_g=r_{l_1}*\cdots*r_{l_k}$,
\item $r_g=\Psi(e_g)+\sum_{w<g}x_{gw}\Psi(e_w)$ where, for a word
$v=[i_1,\ldots,i_k]$, $e_v=e_{i_1}\cdots e_{i_k}$, and
\item $\{r_g|g\in\G\}$ is a basis for $\W$.
\end{enumerate}
\end{thm}

The basis $\{r_g\mid g\in\G\}$ is called the Lyndon basis of $\W$. An immediate
consequence of Proposition \ref{P:bracketingtriangularity} and Theorem
\ref{T:Lyndonbasis} is the following:

\begin{prp}\label{P:LyndonCoveringProperty}\cite[Proposition 24]{lec} Assume $\bt_1,\bt_2\in \Dt^+$,
$\bt_1+\bt_1=\bt\in\Dt^+$, and $l(\bt_1)<l(\bt_2)$. Then, $l(\bt_1)l(\bt_2)\geq
l(\bt)$.
\end{prp}

This gives an inductive algorithm to determine $l(\bt)$ for $\bt\in\Dt^+$ (cf.
\cite[$\S4.3$]{lec}):

For $\af_i\in\Pi\subset\Dt^+$, $l(\af_i)=[i]$. If $\bt$ is not a simple root,
then there exists a factorization $l(\bt)=l_1l_2$ with $l_1,l_2$ Lyndon words.
By Lemma \ref{L:GoodFactors}, $l_1$ and $l_2$ are good, so $l_1=l(\bt_1)$ and
$l_2=l(\bt_2)$ for some $\bt_1,\bt_2\in\Dt^+$ with $\bt_1+\bt_2=\bt$. Assume
that we know $l(\bt_0)$ for all $\bt_0\in\Dt^+$ satisfying
$\height(\bt_0)<\height(\bt)$. Define
\[
C(\bt)=\{\,(\bt_1,\bt_2)\in\Dt^+\times\Dt^+ \mid \bt=\bt_1+\bt_2, \mbox{ and
}l(\bt_1)<l(\bt_2)\,\}.
\]
Then, Proposition \ref{P:LyndonCoveringProperty} implies

\begin{prp}\cite[Proposition 25]{lec} We have
\[
l(\bt)=\min\{\,l(\bt_1)l(\bt_2) \mid (\bt_1,\bt_2)\in C(\bt)\,\}
\]
\end{prp}

\subsection{PBW and Canonical Bases}\label{SS:PBWandCanonical}
The lexicographic ordering on $\GL$ induces a total ordering on $\Dt^+$, which
is \emph{convex}, meaning that if $\bt_1,\bt_2\in\Dt^+$ with $\bt_1<\bt_2$, and
$\bt=\bt_1+\bt_2\in\Dt^+$, then $\bt_1<\bt<\bt_2$ (cf. \cite{ro3,lec}). Indeed,
assume $\bt_1,\bt_2,\bt=\bt_1+\bt_2\in\Dt^+$ and $\bt_1<\bt_2$. Proposition
\ref{P:LyndonCoveringProperty} and \eqref{E:rightfactors} imply that
$l(\bt)\leq l(\bt_1)l(\bt_2)<l(\bt_2)$. If $l(\bt)=l(\bt_1)l(\bt_2)$, then the
definition of Lyndon words implies $l(\bt_1)<l(\bt)$. We are therefore left to
prove that $l(\bt_1)<l(\bt)$ even if $l(\bt)<l(\bt_1)l(\bt_2)$. This can be
checked easily in all cases. We call this ordering a (right) Lyndon ordering on
$\Dt^+$.

Now, \cite[Corollary 27]{lec} becomes

\begin{cor}\label{C:MaxGoodWord} Let $\bt\in\Dt^+$. Then, $l(\bt)$ is the largest good word of weight $\bt$.
\end{cor}

Each convex ordering, $\bt_1<\cdots<\bt_N$, on $\Dt^+$ arises from a unique
decomposition $w_0=s_{i_1}s_{i_2}\cdots s_{i_N}$ of the longest element of the
Weyl group of $\g$ via
\[
\bt_1=\af_{i_1},\;\bt_2=s_{i_1}\af_{i_2},\;\cdots,\bt_N=s_{i_1}\cdots s_{i_{N-1}}\af_{i_N}.
\]
Lusztig associates to this data a PBW basis of $\U_\A$ denoted
\[
E^{(a_1)}(\bt_1)\cdots E^{(a_n)}(\bt_N),\;\;\;(a_1,\ldots,a_N)\in\Z_{\geq0}^N.
\]
Leclerc \cite[$\S4.5$]{lec} describes the image in $\W$ of this basis for the
convex Lyndon ordering. We use the same braid group action as Leclerc and the
results of \cite[$\S4.5,4.6$]{lec} carry over, making changes in the same
manner indicated in the previous section. We describe the relevant facts below.

For $g=l(\bt_1)^{a_1}\cdots l(\bt_k)^{a_k}$, where $\bt_1>\cdots>\bt_k$ and
$a_1,\ldots,a_k\in\Z_{>0}$ set
\[
E_g=\Psi(E^{(a_k)}(\bt_k)\cdots E^{(a_1)}(\bt_1))\in\W_\A
\]
and let $E_g^*\in\W_\A^*$ be the image of $(E^{(a_k)}(\bt_k)\cdots
E^{(a_1)}(\bt_1))^*\in\U_\A^*$. Observe that the order of the factors in the
definition of $E_g$ above are increasing with respect to the Lyndon ordering.
Leclerc shows that if $\bt\in\Dt^+$, then
\begin{align}\label{E:Proportional}
\kp_{l(\bt)}E_{l(\bt)}=r_{l(\bt)},
\end{align}
For some $\kp_{l(\bt)}\in\Q(q)$, \cite[Theorem 28]{lec} (the proof of this
theorem in our case is obtained by reversing all the inequalities and using the
standard factorization as opposed to the costandard factorization). More
generally, let $g=l_1^{a_1}\cdots l_k^{a_k}\in\G$, $l_1>\cdots>l_k\in\GL$. If
$l=l(\bt)$, write $d_l:=d_i$ if $(\bt,\bt)=(\af_i,\af_i)$, and define
\begin{eqnarray}\label{E:kpg}
\kp_g=\prod_{i=1}^k\kp_{l_i}^{a_i}[a_i]_{l_i}!.
\end{eqnarray}
Then, $E_g=\kp_g\sm(r_g)$, where $\sm$ is defined in Propsition
\ref{P:BarInv}, \cite[$\S4.6$]{lec}. Moreover,
\begin{align}\label{E:Normalized Eg}
E_g^*=q^{c_g} (E_{l_m}^*)^{*a_m}*\cdots*(E_{l_1}^*)^{*a_1}
\end{align}
where $c_g=\sum_{i=1}^m{a_i\choose2}d_{l_i}$, \cite[$\S5.5.3$]{lec}.

It is well known that using the bar involution (Proposition \ref{P:BarInv}) we
obtain a canonical basis $\{b_g \mid g\in\G\}$ for $\W_\A$ via the PBW basis
$\{E_g \mid g\in\G\}$, see \cite[Lemma 37]{lec}. It has the form
\begin{align}\label{bg}
b_g=E_g+\sum_{\substack{h\in\G\\h<g}}\chi_{gh}E_h.
\end{align}
The dual canonical basis then has the form
\begin{align}\label{bstarg}
b_g^*=E_g^*+\sum_{\substack{h\in\G\\h>g}}\chi_{gh}^*E_h^*.
\end{align}

As in \cite{lec} we have the following very important theorem:

\begin{thm}\label{T:minbg}\cite[Theorem 40, Corollary 41]{lec}
\begin{enumerate}
\item[(i)] We have $\min(b_g^*)=g$ for all $g\in\G$. Moreover, the coefficient of $g$ in $b_g^*$ is equal to $\kp_g$.
\item[(ii)] For each $l\in\GL$, $E^*_l=b^*_l$.
\end{enumerate}
\end{thm}

To describe the coefficient $\kp_l$ precisely, transport the symmetric
bilinear form \eqref{E:SymmetricForm} to $\W$ via the isomorphism $\Psi$. Let
$g=l(\bt_1)^{a_1}\cdots l(\bt_N)^{a_N}$ and $h=l(\bt_1)^{b_1}\cdots
l(\bt_N)^{b_N}$, where $a_1,\ldots,a_N,b_1,\ldots,b_N\in\Z_{\geq0}$. Then, the
form is given by
\begin{align}\label{E:KashiwaraForm1}
(E_g,E_h)_K=\dt_{gh}\prod_{j=1}^n\frac{(E(\bt_j),E(\bt_j))_K^{a_j}}{\{a_j\}_{(\bt_j,\bt_j)}!}
\end{align}
where, for $\bt=\sum_{i=1}^rc_i\af_i\in\Dt^+$,
\begin{align}\label{E:KashiwaraForm2}
(E(\bt),E(\bt))_K=\frac{\prod_{i=1}^r(1-q^{(\af_i,\af_i)})^{c_i}}{1-q^{(\bt,\bt)}}
\end{align}
and for $a,b\in\Z_{\geq0}$,
\begin{align}\label{E:KashiwaraForm3}
\{a\}_b!=\prod_{j=1}^a\frac{1-q^{jb}}{1-q^b}.
\end{align}
Then, \cite[$\S5.5.2$]{lec},
\begin{align}\label{E:Proportionality}
E_l^*=\frac{(-1)^{\ell(l)-1}\kp_l^{-1}}{q^{N(|l|)}(E_l,E_l)_K}\,r_l,
\end{align}
where $N(|l|)$ is given by \eqref{E:Nfcn}

\subsection{The Anti-Automorphism $\tau$}\label{SS:antiautomorphism of W}
We continue with a fixed ordering, $\leq$, on $I$ and corresponding sets $\G$,
$\L$, and $\GL$ as described in $\S$\ref{SS:LyndonWords}. Define the
\emph{opposite} ordering on $I$ by $$x\preceq y\mbox{ if, and only if, }y\leq
x.$$ Given this opposite ordering, define the corresponding opposite total
ordering on the monomials in $\F$ by
$$[i_1,\ldots,i_k]\prec[j_1,\ldots,j_\ell]\mbox{ if }i_1\prec j_1,\mbox{ or for
some }m, i_m\prec j_m\mbox{ and }i_s=j_s\mbox{ for all }s<m,$$ and $[]\prec[i]$ for
all $i\in I$.

For $f\in\F$, $\max(f)$ is the largest monomial occuring in the expansion of
$f$. Call a monomial $g^\tau=[i_1,\ldots,i_k]\in\F$ an \emph{upper good word}
if $g^\tau=\max(u)$ for some $u\in\U_q$, and we say that it is  \emph{Lyndon
on the left} if it is smaller than all of its proper right factors:
$$[i_1,\ldots,i_k]\prec[i_j,\ldots,i_k]\mbox{ for }j>1.$$
Let $\G^\tau$ denote the set of upper good words, let $\L^\tau$ denote the set
of words that are Lyndon on the left, and $\GL^\tau=\G^\tau\cap\L^\tau$.

Observe that the total ordering on $\GL^\tau$ induces a convex total ordering
on $\Dt^+$ which we call a (left) Lyndon ordering. Also, the bijection
$\Dt^+\to\GL^\tau$ provides a means to compute $l^\tau(\bt)$ for each
$\bt\in\Dt^+$, see \cite[Section 4]{lec}. Finally, given $l^\tau\in\L^\tau$,
define its costandard factorization to be $l^\tau=l_1^\tau l_2^\tau$, where
$l_1^\tau$ is the maximal proper word which is Lyndon on the left. Note that
$l_2^\tau$ is also Lyndon on the left. Using the data above we may define a
Lyndon basis $\{r_{g^\tau}|g^\tau\in\G^\tau\}$, dual PBW basis
$\{E^*_{g^\tau}\,|\,g^\tau\in\G^{\tau}\}$ and dual canonical basis
$\{b^*_{g^\tau}\,|\,g^{\tau}\in\G^\tau\}$ exactly as in \cite[Sections 4-5]{lec}.

The next lemma gives the precise connection between the combinatorics appearing
here and those developed in \cite{lec}:

\begin{lem}\label{L:antiautomorphism of W} Under the anti-automorphism $\tau:\F\to\F$,
\begin{enumerate}
\item $\tau(\W)=\W$;
\item $\tau(\G)=\G^\tau$ and $\tau(\L)=\L^\tau$;
\item $\tau(E^*_g)=E^*_{\tau(g)}$;
\item $\tau(b^*_g)=b^*_{\tau(g)}$.
\end{enumerate}
\end{lem}

\begin{pff}
Property (1) is immediate from Proposition \ref{P:BarInv}, and property (2) is
clear from the definitions.

We now turn to property (3). Observe that if $g=l_1\cdots l_k$, then
$\tau(g)=\tau(l_k)\cdots\tau(l_1)$. Therefore, by equation \eqref{E:Normalized
Eg}, it is enough to show that $\tau(E^*_l)=E^*_{\tau(l)}$ for all $l\in\GL$.
We prove this by induction on the degree of $l$ in the principal grading on
$\F$. The base case is clear since $E^*_{[i]}=r_{[i]}=[i]$.

For the inductive step, assume we have shown that $\tau(r_{l_0})=r_{\tau(l_0)}$
and $\tau(E^*_{l_0})=E^*_{\tau(l_0)}$ for all $l_0$ of degree less than the
degree of $l$. Let $l=l_1l_2$ be the standard factorization of $l$. Then, by
(2), $\tau(l)=\tau(l_2)\tau(l_1)$ is the costandard factorization of $\tau(l)$.
Then, it follows from \eqref{E:Lyndonbracketing} and the relevant definitions
that
\begin{align*}
\tau(r_l)&=\tau(r_{l_1}*r_{l_2}-q^{(|l_1|,|l_2|)}r_{l_2}*r_{l_1})\\
    &=\tau(r_{l_2})*\tau(r_{l_1})-q^{(|l_1|,|l_2|)}\tau(r_{l_1})*\tau(r_{l_2})\\
    &=r_{\tau(l_2)}*r_{\tau(l_1)}-q^{(|l_2|,|l_1|)}r_{\tau(l_1)}*r_{\tau(l_2)}\\
    &=r_{\tau(l)}.
\end{align*}
It now follows that $\tau(E^*_l)=E^*_{\tau(l)}$ by applying $\tau$ to
equation \eqref{E:Proportionality} and observing that equations
\eqref{E:Nfcn} and \eqref{E:KashiwaraForm1}-\eqref{E:KashiwaraForm3} imply that the coefficient on the
right-hand-side of \eqref{E:Proportionality} depend only on $|l|\in\curlyQ^+$.

Finally, property (4) for follows by applying $\tau$ to equation \eqref{bstarg} and uniqueness.
\end{pff}

From now on, we will write $g^\tau=\tau(g)$.

\section{Quiver Hecke Algebras}\label{S:QHA}

In this section, we give a presentation of the quiver Hecke algebras following
the notation of \cite{klram2}. Throughout, we work over an arbitrary ground
field $\bbF$.

\subsection{Quivers with Compatible Automorphism}\label{SS:QuivwithAut}
Let $\widetilde{\Gm}$ be a graph. We construct a Dynkin diagram $\Gm$ by giving
$\widetilde{\Gm}$ the structure of a \emph{graph with compatible automorphism}
in the sense of \cite[$\S12,14$]{lu}. To define the quiver Hecke algebra, we
will use the notion of a \emph{quiver with compatible automorphism} as
described in \cite[$\S3.2.4$]{rq}.

Let $\widetilde{I}$ be the labelling set for $\widetilde{\Gm}$, and
$\widetilde{H}$ be the (multi)set of edges. An automorphism
$a:\widetilde{\Gm}\to\widetilde{\Gm}$ is said to be \emph{compatible} with
$\widetilde{\Gm}$ if, whenever $(i,j)\in \widetilde{H}$ is an edge, $i$ is not
in the orbit of $j$ under $a$.

Fix a compatible automorphism $a:\widetilde{\Gm}\to\widetilde{\Gm}$, and set
$I$ to be a set of representatives of the obits of $\widetilde{I}$ under $a$
and, for each $i\in I$, let $\af_i\in\widetilde{I}/a$ be the corresponding
orbit. For $i,j\in I$, $i\neq j$ define $(\af_i,\af_i)=2|\af_i|$ and let
$$(\af_i,\af_j)=-|\{(i',j')\in\tilde{H}\,|\,i'\in\af_i,j'\in\af_j\}|.$$ For all
$i,j\in I$, let $a_{ij}=2(\af_i,\af_j)/(\af_i,\af_i)$. Then, \cite[Proposition
14.1.2]{lu} $A=(a_{ij})_{i,j\in I}$ is a Cartan matrix and every Cartan matrix
arises in this way. Let $\Gm$ be the Dynkin diagram corresponding to $A$.
Moreover, the pairing $(\af_i,\af_j)$ defined above agrees with the pairing on
$\curlyQ$ in $\S$\ref{SS:Root Data}.

Assume further that $\widetilde{\Gm}$ is a quiver. That is, we have a pair of
maps $s:\widetilde{H}\to\widetilde{I}$ and $t:\widetilde{H}\to\widetilde{I}$
(the source and the target). We say that $a$ is a compatible automorphism if
$s(a(h))=a(s(h))$ and $t(a(h))=a(t(h))$ for all $h\in\widetilde{H}$. Set
$$d_{ij}=|\{h\in\widetilde{H}\,|\,s(h)\in\af_i\mbox{ and }t(h)\in\af_j\}/a|$$
and let $m(i,j)=\mathrm{lcm}\{(\af_i,\af_i),(\af_j,\af_j)\}$. As noted in
\cite{rq},
\begin{align}\label{E:dij+dji}
d_{ij}+d_{ji}=-2(\af_i,\af_j)/m(i,j).
\end{align}

This data defines a matrix $Q=(Q_{ij}(u,v))_{i,j\in I}$, where each
$Q_{ij}(u,v)\in\bbF[u,v]$. The polynomial entries in $Q$ are defined by
$Q_{ii}(u,v)=0$, and for $i\neq j$,
\begin{align}\label{E:Q polynomials}
Q_{ij}(u,v)=(-1)^{d_{ij}}(u^{m(i,j)/(\af_i,\af_i)}-v^{m(i,j)/(\af_j,\af_j)})^{-2(\af_i,\af_j)/m(i,j)}
\end{align}

Specialize now to the case where $\Gm$ is of finite type. Then, as explained in
\cite[$\S3.1$]{klram2}, the polynomials $Q_{ij}(u,v)$ ($i\neq j$) are
completely determined by the Cartan matrix and a partial ordering on $I$ such
that $i\to j$ or $j\to i$ if $a_{ij}\neq0$. In this case,
\begin{align}\label{E:finite type Q polynomials}
Q_{ij}(u,v)=\begin{cases}0&\mbox{if }i=j;\\
    1&\mbox{if }a_{ij}=0;\\u^{-a_{ij}}-v^{-a_{ji}}&\mbox{if }a_{ij}<0\mbox{ and }i\to j;\\
    v^{-a_{ji}}-u^{-a_{ij}}&\mbox{if }a_{ij}<0\mbox{ and }j\to i.\end{cases}
\end{align}

\subsection{Generators and Relations}\label{SS:QHAGenRel}

Assume from now on that $\g$ is as in $\S$\ref{SS:Root Data}. Define the quiver
Hecke algebra
\[
H(\Gm)=\bigoplus_{\nu\in\curlyQ^+}H(\Gm;\nu),
\]
where $H(\Gm;\nu)$ is the unital $\bbF$-algebra, with identity $1_\nu$, given by generators and relations as
described below.

Assume that $\height(\nu)=d$. The set of generators are
\[
\{e(\ui)|\ui\in I^\nu\}\cup\{y_1,\ldots, y_d\}\cup\{\phi_1,\ldots,\phi_{d-1}\}.
\]
We refer to the $e(\ui)$ as \emph{idempotents}, the $y_r$ as \emph{Jucys-Murphy
elements}, and the $\phi_r$ as \emph{intertwining elements}. Indeed, these
generators are subject to the following relations for all $\ui,\uj\in I^\nu$
and all admissible $r,s$:
\begin{align}
\label{E:Relations1}e(\ui)e(\uj)&=\dt_{\ui,\uj}e(\ui);\\
\label{E:Relations2}\sum_{\ui\in I^\nu}e(\ui)&=1_\nu;\\
\label{E:Relations3}y_re(\ui)&=e(\ui)y_r;\\
\label{E:Relations4}\phi_re(\ui)&=e(s_r\cdot\ui)\phi_r;\\
\label{E:Relations5}y_ry_s&=y_sy_r;\\
\label{E:Relations6}\phi_ry_s&=y_s\phi_r&\mbox{if }s\neq r,r+1;\\
\label{E:Relations7}\phi_r\phi_s&=\phi_s\phi_r&\mbox{if }|s-r|>1;\\
\label{E:Relations8}\phi_ry_{r+1}e(\ui)&=\begin{cases}(y_r\phi_r+1)e(\ui)\\y_r\phi_re(\ui)\end{cases}&\begin{array}{r}i_r=i_{r+1},\\i_r\neq
                        i_{r+1};\end{array}\\
\label{E:Relations9}y_{r+1}\phi_re(\ui)&=\begin{cases}(\phi_ry_r+1)e(\ui)\\\phi_ry_re(\ui)\end{cases}&\begin{array}{r}i_r=i_{r+1},\\i_r\neq
                        i_{r+1}.\end{array}
\end{align}
Additionally, the intertwining elements satisfy the quadratic relations
\begin{align}\label{E:Quadratic Relations}
\phi_r^2e(\ui)=Q_{i_r,i_{r+1}}(y_r,y_{r+1})e(\ui)
\end{align}
for all $0\leq r<d$, and the braid-like relations
\begin{align}\label{E:Braid Relations}
(\phi_r\phi_{r+1}\phi_r-\phi_{r+1}&\phi_r\phi_{r+1})e(\ui)\\
\nonumber&=\begin{cases}\left(\frac{Q_{i_r,i_{r+1}}(y_{r+2},y_{r+1})-Q_{i_r,i_{r+1}}(y_r,y_{r+1})}{y_{r+2}-y_r} \right)e(\ui)&\mbox{if }i_r=i_{r+2},\\
    0&\mbox{otherwise.}\end{cases}
\end{align}
Finally, this algebra is graded via
\begin{align}\label{E:Grading}
\deg e(\ui)=0,\;\;\;\deg
y_re(\ui)=(\af_{i_r},\af_{i_r}),\andeqn\deg\phi_re(\ui)=-(\af_{i_r},\af_{i_{r+1}}).
\end{align}

\subsection{Basis Theorem}\label{SS:BasisThm} Let $\nu\in\curlyQ^+$ with $\height(\nu)=d$.
Given $w\in S_d$, fix a reduced decomposition $w=s_{k_1}\cdots s_{k_t}$ for $w$
and define
\[
\phi_w=\phi_{k_1}\cdots\phi_{k_t}.
\]
Relations \eqref{E:Quadratic Relations} and \eqref{E:Braid Relations} imply
that, in general, $\phi_w$ depends on the choice of reduced decomposition.



Finally, we have

\begin{thm}\label{T:Basis Theorem}\cite[Theorem 2.5]{khl1}\cite[Theorem 3.7]{rq} The set
\[
\{\,\phi_wy_1^{m_1}\cdots y_d^{m_d}e(\ui)\,|\,w\in
S_d,\,m_1,\ldots,m_d\in\Z_{\geq0},\,\ui\in I^\nu\,\}
\]
forms an $\bbF$-basis for $H(\Gm;\nu)$.
\end{thm}


\subsection{An Automorphism and Anti-Automorphism of $H(\Gm;\nu)$}\label{SS:Automorphism}
Let $\nu\in\curlyQ^+$, $\height(\nu)=d$. As observed in \cite[$\S$2.1]{khl1}, we have the following

\begin{prp}\label{P:AutomorphismTau} There is a unique $\bbF$-linear automorphism
$\tau:H(\Gm;\nu)\to H(\Gm;\nu)$ given by
$\tau(e(i_1,\ldots,i_d))=e(i_d,\ldots,i_1)$, $\tau(y_r)=y_{d-r+1}$, and
$\tau(\phi_r)=-\phi_{d-r}$.
\end{prp}

and

\begin{prp}\label{P:Anti-AutomorphismPsi} There is a unique $\bbF$-linear anti-automorphism $\psi:H(\Gm;\nu)\to H(\Gm;\nu)$ defined by $\psi(e(\ui))=e(\ui)$, $\psi(y_r)=y_r$, and $\psi(\phi_r)=\phi_r$ for all $\ui\in I^\nu$ and admissible $r$.
\end{prp}

\subsection{Modules and Graded Characters}\label{SS:ModulesandChar} Given a finite dimensional
$\Z$-graded vector space $V=\bigoplus_{k\in\Z}V[k]$, define the \emph{graded
dimension} of $V$ to be
\[
\gdim V=\sum_{k\in\Z}(\dim V[k])q^k\in\Z_{\geq0}[q,q^{-1}].
\]
Let $V\!\{s\}$ denote the vector space obtained from $V$ by shifting the grading by $s$. That is,
\[
\gdim V\!\{s\}=q^s\gdim V.
\]

The algebra $H(\Gm;\nu)$ is $\Z$-graded by \eqref{E:Grading}. Let $\Rep(\Gm;\nu)$ denote the category of all finite dimensional graded
$H(\Gm;\nu)$-modules. Let $M$ be in $\Rep(\Gm;\nu)$. For each $\ui\in I^\nu$,
define the generalized $\ui$-eigenspace by $M_\ui:=e(\ui)M$. We have the
decomposition
\[
M=\bigoplus_{\ui\in I^\nu}M_\ui.
\]
Moreover, by \eqref{E:Relations4}, $\phi_rM_\ui=M_{s_r\cdot\ui}$. Finally, note
that since the elements $y_re(\ui)$ have positive degree, they act nilpotently
on all objects in $\Rep(\Gm;\nu)$.

Morphisms are \emph{degree 0} $H(\Gm;\nu)$-homomorphisms. That is,
for each $M,N\in\Rep(\Gm;\nu)$, $\hom_\nu(M,N)$ denotes the set of \emph{degree 0} homomorphisms.

Let $K(\Gm;\nu)=K(\Rep(\Gm;\nu))$ be the Grothendieck group of the category
$\Rep(\Gm;\nu)$, and
\[
K(\Rep(\Gm)):=K(\Gm)=\bigoplus_{\nu\in\curlyQ^+}K(\Gm;\nu).
\]
This is a free $\Z[q,q^{-1}]$-module with basis given by isomorphism classes of
simple $H(\Gm)$-modules. Note that since morphisms have degree 0, $L\ncong
L\!\{s\}$ for any simple module $L\in\Rep(\Gm;\nu)$ and any $s\neq 0$. We
write $[M]\in K(\Gm;\nu)$ for the image of $M\in\Rep(\Gm;\nu)$ in the
Grothendieck group. Finally, observe that $q^s[M]=[M\!\{s\}]$. Define the
\emph{formal character} $\ch:\Rep(\Gm;\nu)\rightarrow\F$ by
\[
\ch M=\sum_{\ui\in I^\nu}(\gdim M_\ui)\cdot[\ui].
\]

\begin{thm}\label{T:Independence of Characters}\cite[Theorem 3.17]{khl1}
The character map induces an injective $\Q(q)$-linear map
\[
\ch:K(\Gm;\nu)\rightarrow\F.
\]
\end{thm}

Now, let $\nu,\nu'\in\curlyQ^+$ and let $H(\Gm;\nu,\nu'):=H(\Gm;\nu)\otimes
H(\Gm;\nu')$. Given $\ui=(i_1,\ldots,i_{d})\in I^\nu$ and
$\uj=(j_1,\ldots,j_{d'})\in I^{\nu'}$, let
$\ui\uj=(i_1,\ldots,i_d,j_1,\ldots,j_{d'})$. Then, there exists an embedding
\begin{align}\label{E:iota}
\io_{\nu,\nu'}:H(\Gm;\nu,\nu')\rightarrow H(\Gm;\nu+\nu')
\end{align}
given by $\io_{\nu,\nu'}(e(\ui)\otimes e(\uj))=e(\ui\uj)$ and, for appropriate
$r$ and $s$, and for $a$ and $b$ among the symbols $y$ or $\phi$,
\[
\io_{\nu,\nu'}(a_r\otimes b_s)=a_rb_{s+d}.
\]

If $M\in\Rep(\Gm;\nu)$ and $N\in \Rep(\Gm;\nu')$, let $M\boxtimes N\in\Rep(\Gm;\nu)\otimes \Rep(\Gm;\nu')$ denote the \emph{outer} tensor product of $M$
and $N$. We have

\begin{prp}\cite[Proposition 2.16]{khl1} We have
$\io_{\nu,\nu'}(1_\nu\otimes1_{\nu'})H(\Gm;\nu+\nu')$ is a free graded left $H(\Gm;\nu,\nu')$-module.
\end{prp}

Therefore, we may define the exact functors
\begin{align}\label{E:Res}
\res^{\nu+\nu'}_{\nu,\nu'}:H(\Gm;\nu+\nu')\to H(\Gm;\nu,\nu')
\end{align}
by $\res^{\nu+\nu'}_{\nu,\nu'}M=\io_{\nu,\nu'}(1_\nu\otimes 1_{\nu'})M$, and
\begin{align}\label{E:Ind}
\ind_{\nu,\nu'}^{\nu+\nu'}:H(\Gm;\nu,\nu')\rightarrow
H(\Gm;\nu+\nu'),
\end{align}
by
\[
\ind_{\nu,\nu'}^{\nu+\nu'}M\boxtimes
N=H(\Gm;\nu+\nu')\otimes_{H(\Gm;\nu,\nu')}M\boxtimes N.
\]
We have

\begin{lem}\label{P:Shuffle Lemma}\cite[Lemma 2.20]{khl1} Assume that
$M\in\Rep(\Gm;\nu)$, $N\in \Rep(\Gm;\nu')$,
\[
\ch M=\sum_{\ui\in I^\nu}m_\ui[\ui]\andeqn\ch N=\sum_{i\in
I^{\nu'}}n_{\uj}[\uj].
\]
Then,
\[
\ch\ind_{\nu,\nu'}^{\nu+\nu'}M\boxtimes N=\sum_{\ui\in I^\nu,i\in
I^{\nu'}}m_\ui n_\uj[\uj]*[\ui]
\]
where $[\uj]*[\ui]$ is the shuffle product given by \eqref{E:Shuffle}.
\end{lem}

\begin{rmk}\label{R:different shuffles} Observe that the order of the segments in the shuffle lemma is reversed.
This is a consequence of the definition \eqref{E:inductiveqshuffle} and is so
that the terms in the character formula coming from $1\otimes (M\boxtimes N)$
are not shifted in degree. Note that this is slightly different than the
shuffle product in \cite{klram2}. The products are related by the formula
\[
x\circ y=y*x
\]
for $x,y\in\W$.
\end{rmk}

Let $\Proj(\Gm)$ (resp. $\Proj(\Gm;\nu)$) denote the category
of finitely generated, graded, projective (left) $H(\Gm)$-modules (resp.
$H(\Gm;\nu)$-modules). Let $K_0(\Gm)$ (resp. $K_0(\Gm;\nu)$) denote Grothendieck group of $\Proj(\Gm)$ (resp. $\Proj(\Gm;\nu)$), and set $K_0(\Gm)_{\Q(q)}=K_0(\Gm)\otimes_\A\Q(q)$. Also, define the analogous functors $\ind_{\nu,\nu'}^{\nu+\nu'}$ and $\res_{\nu,\nu'}^{\nu+\nu'}$ to those for $\Rep(\Gm)$.

Given a left $H(\Gm;\nu)$-module $M$, let $M^\psi$ be the right $H(\Gm;\nu)$-module given by $mx=\psi(x)m$ for all $x\in H(\Gm;\nu)$ and $m\in M$. Define the \emph{Kashiwara-Khovanov-Lauda pairing} (KKL), $(\cdot,\cdot)_{KKL}:K_0(\Gm;\nu)\times K(\Gm;\nu)\to\A$ by \begin{align}\label{E:RenormalizedKLform}
([P],[M])_{KKL}=\prod_{i\in I}(1-q_i^2)^{c_i}\,\gdim(P^\psi\otimes M),
\end{align}
if $\nu=\sum_{i\in I}c_i\af_i$. This form is evidently related to the \emph{Lusztig-Khovanov-Lauda pairing} (LKL), $(\cdot,\cdot)_{LKL}$, appearing in \cite[(2.43),(2.44)]{khl1} by the formula
\begin{align}\label{E:RenormalizedKLform2}
([P],[M])_{KKL}=\prod_{i\in I}{(1-q_i^2)^{c_i}}\,([P],[M])_{LKL},
\end{align}
see Remark \ref{R:LusztigForm}. Define the map $\omega:K(\Gm)\to K_0(\Gm)_{\Q(q)}$,
\begin{align}\label{E:omega}
\omega([M])=\sum_{[P]\in\mathcal{B}}([P],[M])_{KKL}[P],
\end{align}
where the sum is over a basis $\mathcal{B}$ of $K_0(\Gm)$, and $M\in\Rep(\Gm)$.

\begin{exa}\label{Exa:KKLform}
Let $\va_{\af_i}$ denote the unique irreducible $H(\Gm;\af_i)$-module concentrated in degree 0. It is one dimensional with the action of $H(\Gm;\af_i)$ given by $e(j)\va_{\af_i}=\dt_{ij}\va_{\af_i}$, $y_1\va_{\af_i}=\{0\}$. Let $\P_{\af_i}$ denote its projective cover. Then,
$$([\P_{\af_i}],[\va_{\af_i}])_{KKL}=(1-q_i^2).$$
In particular, under the identification of $K(\Gm)$ with $K^*_0(\Gm)$, we have $\omega([\va_{\af_i}])=[\P_{\af_i}]-[\P_{\af_i}\{2d_i\}]$. That is, $[\va_{\af_i}]$ is mapped by $\omega$ to its projective resolution $$\xymatrix{0\ar[r]&\P_{\af_i}\{2d_i\}\ar[r]&\P_{\af_i}\ar[r]&\va_{\af_i}\ar[r]&0}.$$

More generally, using \cite[Lemma 3.2]{klram2}, we deduce that if $\L\in\Rep(\Gm;\nu)$ is a simple module satisfying $\sm(\ch\L)=\ch\L$, $\P_{\L}\in\Proj(\Gm;\nu)$ is its projective cover, and $\nu=\sum_ic_i\af_i$, then
$$([\P_\L],[\L])_{KKL}=\prod_{i\in I}(1-q_i^2)^c_i,$$ so $\omega([\L])=\prod_i(1-q_i^2)^c_i[\P_\L]$. On the other hand, if $\L,\L'\in\Rep(\Gm)$ are two simple modules as above, $$([\P_{\L}],[\L'])_{LKL}=\dt_{\L,\L'},$$ so identifying $K(\Gm)$ with the dual lattice to $K_0(\Gm)$ using the Lusztig-Khovanov-Lauda pairing does not contain any representation theoretic information.
\end{exa}

We identify $K(\Gm)$ with its image under $\omega$. The following lemma shows that this image is the dual lattice $$K^*_0(\Gm)=\{X\in K_0(\Gm)_{\Q(q)}\,|\,(Y,X)_{LKL}\in\A\mbox{ for all }Y\in K_0(\Gm)\},$$ where 
\begin{align}\label{E:KLform}
(\cdot,\cdot)_{LKL}:K_0(\Gm)_{\Q(q)}\times K_0(\Gm)_{\Q(q)}\to\Q(q)
\end{align}
is the \emph{Lusztig-Khovanov-Lauda} bilinear form, given by $([P],[Q])_{LKL}=\gdim(P^\psi\otimes_{H(\Gm;\nu)}Q)$ for $P,Q\in\Proj(\Gm;\nu)$, cf. \cite[(2.45),(2.46),(2.47)]{khl1}.

\begin{lem}\label{L:IdentificationofLattices}
Under the identification above, the simple modules are dual to their projective covers with respect to the Lusztig-Khovanov-Lauda bilinear form. In particular, the map $X\mapsto (\omega(X),?)_{LKL}$ identifies the dual lattice $K^*(\Gm;\nu)$ with the dual space $\Hom_{\A}(K_0(\Gm;\nu),\A)$.
\end{lem}

\begin{pff}
Let $\nu=\sum_ic_i\af_i\in\curlyQ^+$. Assume that $\{\L_a|a\in A\}$ is a basis for $K(\Gm;\nu)$ for some indexing set $A$, and let $\mathcal{B}$ in \eqref{E:omega} be the basis for $K_0(\Gm;\nu)$ consisting of the projective covers $\P_a$ of $\L_a$, $a\in A$. Then, by the definitions $([\P_a],[\L_b])_{KKL}=\dt_{ab}\prod_i(1-q_i^2)^{c_i}$. Therefore, $\omega([\L_a])=\prod_i(1-q_i^2)^{c_i}[\P_a]$. Also, $([\P_b],[\P_a])_{LKL}=\dt_{ba}(1-q_i^2)^{-c_i}$. Hence,
$$(\omega([\L_b]),[\P_a])_{LKL}=(([\P_b],[\L_b])_{KKL}[\P_b],[\P_a])_{LKL}=([\P_b],[\L_b])_{KKL}([\P_b],[\P_a])_{LKL}=\dt_{ba}.$$
\end{pff}

Identifying $\U_{\A,\nu}^*$ with $\Hom_\A(\U_{\A,\nu},\A)$ using Kashiwara's bilinear form, we obtain the following result which is dual to the main results in \cite{khl1,khl2}:

\begin{thm}\label{T:Isomorphism of Bialgebras}\cite[Theorem 1.1]{khl1},\cite[Theorem 8]{khl2}
In the notation of $\S$\ref{SS:Root Data}-\ref{SS:Embedding}, there is an
isomorphism of $\curlyQ^+$-graded twisted bialgebras $$\gm^*:K(\Gm)\to\U_\A^*.$$
\end{thm}

Define multiplication $\circ:K(\Gm;\nu)\otimes K(\Gm;\nu')\to
K(\Gm;\nu+\nu')$ by
$$[M]*[N]=[\ind_{\nu',\nu}^{\nu+\nu'}N\boxtimes M].$$ Define multiplication on $\Proj(\Gm)$ by $[P]\cdot[Q]=[\ind_{\nu,\nu'}^{\nu+\nu'}P\boxtimes Q]$. In light of Remark \ref{R:different shuffles}, we have the following slight modification of \cite[Lemma 3.5]{klram2}:

\begin{lem}\label{L:Coproduct}\cite[Lemma 3.5]{klram2} For $P\in\Proj(\Gm;\nu+\nu')$, $M\in\Rep(\Gm;\nu)$ and $N\in\Rep(\Gm;\nu')$, $$([P],[M]*[N])_{KKL}=([\res^{\nu+\nu'}_{\nu',\nu}P], [N]\otimes[M])_{KKL}.$$
For $P\in\Proj(\Gm;\nu)$, $Q\in\Proj(\Gm;\nu')$, and $M\in\Rep(\Gm;\nu+\nu')$,
$$([P]\cdot[Q],[M])_{KKL}=([P]\otimes[Q],[\res^{\nu+\nu'}_{\nu,\nu'}M])_{KKL}.$$
\end{lem}

\begin{rmk}\label{R:Opp order}
We note that using \eqref{E:RenormalizedKLform2} does not affect the lemma above, since the renormalization factor on both sides of the equations above is the same.
\end{rmk}

Observe that the order of $\nu$ and $\nu'$ in the first equation of Lemma \ref{L:Coproduct} have been reversed, but not in the second. This implies that $\gm^*\circ \mathrm{mult}=\mathrm{mult}\circ(\gm^*\otimes\gm^*)\circ\mathrm{flip}$, where $\mathrm{mult}$ denotes the appropriate multiplication map, and $\mathrm{flip}:K(\Gm)\otimes K(\Gm)\to K(\Gm)\otimes K(\Gm)$ is the map $\mathrm{flip}([M]\otimes[N])=[N]\otimes[M]$. In particular, we have the following property of $\gm^*$ as proved in \cite{klram2}.

\begin{thm}\label{T:Prod in Iso of Bialgebras}\cite[Theorem 4.4(5)]{klram2}
For $[M],[N]\in K(\Gm)$, $$\gm^*([M]*[N])=\gm^*([M])\gm^*([N]).$$
\end{thm}

We also record the following, which was proved in \cite{klram2}.

\begin{thm}\label{T:commutative diagram}\cite[Theorem 4.4(3)]{klram2} The following diagram commutes:
\[
\xymatrix{K(\Gm)\ar[dr]_\ch\ar[rr]^{\gm^*}&&\U_\A^*\ar[dl]^{\Psi}\\&\W_\A^*&}
\]
\end{thm}

\begin{pff}
It is more convenient to show that $\ch\circ(\gm^*)^{-1}=\Psi$. To this end,
assume that $u\in\U_{\A,\nu}^*$. Then, $u$ may be written as $$u=\sum
n_{i_1,\ldots,i_d}e_{i_1}\cdots e_{i_d},$$ where the sum is over all
$(i_1,\ldots,i_d)\in I^\nu$.

Now, let $\va_{\af_i}\in\Rep(\Gm;\af_i)$ be the unique irreducible
representation, see Example \ref{Exa:KKLform}. It is clear from Theorem \ref{T:Isomorphism of Bialgebras} that
$\gm^*([\va_{\af_i}])=e_i$. Therefore,
\begin{align*}
\ch\circ(\gm^*)^{-1}(u)&=\ch\circ(\gm^*)^{-1}\left(\sum n_{i_1,\ldots,i_d}e_{i_1}\cdots e_{i_d}\right)\\
    &=\ch\left(\sum n_{i_1,\ldots,i_d}[\va_{\af_{i_1}}]*\cdots*[\va_{\af_{i_d}}]\right)\\
    &=\ch\left(\sum n_{i_1,\ldots,i_d}[\ind_{\af_{i_1},\ldots,\af_{i_d}}^\nu \va_{\af_{i_d}}\boxtimes\cdots\boxtimes\va_{\af_{i_1}}]\right)\\
    &=\sum n_{i_1,\ldots,i_d}[i_1]*\cdots*[i_d]\\
    &=\Psi(u).
\end{align*}
\end{pff}







\begin{rmk}
We point out that Kleshchev and Ram prove several other important properties of the isomorphism $\gm^*$ in \cite[Theorem 4.4]{klram2}. However, as we do not use these properties, we refer the reader to their paper for the details.
\end{rmk}

\section{Standard Representations and their Simple Quotients}\label{S:StdRepsandSimples}

\subsection{Cuspidal Representations}\label{SS:CuspidalReps}
Following Kleshchev and Ram, we call a monomial $f\in\F$ a \emph{weight} of
$M\in\Rep(\Gm)$ if $M_{\ui_f}\neq0$, where $\ui_f\in I^\infty$ is the reading
of the work $f$. That is, $f=[\ui_f]$. Since the set of words in $\F$ is
totally ordered, it makes sense to speak of the \emph{lowest weight} of a
module.

Fix a (right) Lyndon ordering on $\Dt^+$. Continuing with the terminology of
Kleshchev and Ram, we call an irreducible module \emph{cuspidal} if it has
lowest weight $l(\bt)\in\GL$ for some $\bt\in\Dt^+$.

\begin{thm}\label{T:CuspidalsExist}
For the (right) Lyndon order on $\Dt^+$ used in Section
\ref{S:IdentGoodLyndon}, cuspidal representations exist in all finite types.
Moreover, for each $l\in\GL$, $\ch\va_l=b^*_l$.
\end{thm}

\begin{pff}
For types $ABCDF$, the representations are constructed explicitly in Section
\ref{A:SegReps}. We deduce the $E_8$ case from \cite[Lemma 3.3, Theorems
3.6,3.10]{klram}, since the corresponding Lyndon words are homogeneous.
Finally, the $G_2$ case follows easily from the construction in \cite{klram2}
since the characters are identical.
\end{pff}

\subsection{Standard Representations and Unique Irreducible Quotients}\label{SS:StdReps}
We continue to use the ordering from Section \ref{S:IdentGoodLyndon}. Given
$g\in\G$, $g=l(\bt_1)\cdots l(\bt_k)$, with $\bt_1\geq\cdots\geq\bt_k$ define
\[
\M(g)=(\ind_{\bt_1,\ldots,\bt_k}^{\bt_1+\cdots+\bt_k}
\va_{\bt_1}\boxtimes\cdots\boxtimes\va_{\bt_k})\{c_g\}.
\]
The following is a consequence of Lemma \ref{P:Shuffle Lemma},
\eqref{E:Normalized Eg} and the definition.

\begin{prp}\label{P:Std Module Character}
For each $g\in\G$,
\[
\ch\M(g)=E_g^*.
\]
In particular, $\gdim \M(g)_{\ui_g}=\kp_g$.
\end{prp}

The next theorem now follows from the previous proposition using Theorem
\ref{T:commutative diagram}.

\begin{thm}\label{T:StdModBasis}
The set $$\{[\M(g)]\,|\,g\in\G\}$$ forms a basis for $K(\Gm)$.
\end{thm}

The following crucial lemma is proved in \cite{klram2}.

\begin{lem}\cite[Lemma 6.6]{klram2}\label{L:rectangular case}
Let $g=l^k$ for some $l=l(\bt)\in\GL_d$, then $\M(g)$ is irreducible.
\end{lem}

The above lemma, together with a Frobenius reciprocity argument yields the main
result of \cite{klram2}:

\begin{thm}\cite[Theorem 7.2(i)]{klram2}\label{T:Unique Simple Quotient}
Let $g\in\G_d$. Then $\M(g)$ has a unique maximal submodule $\mathcal{R}(g)$
and unique simple quotient $\L(g)$.
\end{thm}

As noted in \cite{klram2}, Khovanov and Lauda prove that for every simple
module $L$, there is a unique grading shift such that $\sm(\ch L\{s\})=\ch
L\{s\}$, \cite[$\S3.2$]{khl1}. Therefore, by Theorems \ref{T:CuspidalsExist}
and \ref{T:minbg}, and \cite[Proposition 32]{lec},

\begin{thm}\cite[Theorem 7.2(iii)]{klram2}
We have $\sm(\ch\L(g))=\ch\L(g)$.
\end{thm}

Finally, we have

\begin{thm}\cite[Theorem 7.2(iv)]{klram2}\label{T:IrrepBasis}
The set $$\{[\L(g)]\,|\,g\in\G\}$$ forms a basis for $K(\Gm)$.
\end{thm}


\subsection{Twisting by the Automorphism $\tau$}\label{SS:Twisting}
Finally, we close by relating the representation theory coming from the (right)
Lyndon orderings on $\Dt^+$ to the (left) Lyndon orderings that appear in
\cite{klram2}. To this end, given $M\in\Rep(\Gm)$, let $M^\tau$ be the module
obtained by twisting by the automorphism $\tau$, cf. Proposition
\ref{P:AutomorphismTau}. That is, $M^\tau=M$ as graded vector spaces with
$x\cdot m=\tau(x)m$ for all $m\in M^\tau$.

Recall the opposite ordering and related notation developed in
$\S$\ref{SS:antiautomorphism of W}. We have the following:

\begin{thm}\label{T:TwistingSimples} Let $g\in\G$. Then, $\L(g)^\tau=\L(g^\tau)$.
\end{thm}

\begin{pff}
First, it is immediate by character considerations that the cuspidal
representations satisfy $\va_l^\tau=\va_{l^\tau}$, see Lemma
\ref{L:antiautomorphism of W}(4). Therefore, it follows that
$\M(g)^\tau=\M(g^\tau)$ for all $g\in\G$. The result now follows since
$\mathcal{R}$ is a submodule of $\M(g)$ if, and only if, $\mathcal{R}^\tau$ is
a submodule of $\M(g^\tau)$.
\end{pff}

\section{Identification of Good Lyndon Words and Associated Root Vectors}\label{S:IdentGoodLyndon}

We now give explicit descriptions of the good Lyndon words and associated root
vectors for $\g$ of classical type and type $F_4$. In type $E_8$ we determine
the good Lyndon words. Throughout, we write $b^*[\ui]:=b^*_{[\ui]}$ for good
Lyndon words $l=[\ui]$.

\subsection{Classical Type}\label{SS:Classicaltype} We now
specialize to the case where $\g$ is of classical type. Fix a labeling of the
simple roots as in Table \ref{Tbl:Dynkin}.

\begin{table}[ht]\label{Tbl:Dynkin}\caption{Labelling of Simple Roots}
\begin{center}
\begin{tabular}{|c|c|c|}\hline
Type&Diagram&Positive Roots\\\hline\hline
$A_r$&
\begin{picture}(76,15)
\multiput(5,5)(21,0){4}{\circle{2}} \multiput(6,5)(42,0){2}{\line(1,0){19}}
\put(30,5){$\ldots$}
\put(4,-2){\tiny$0$}\put(25,-2){\tiny$1$}\put(42,-2){\tiny$r$-$2$}\put(63,-2){\tiny$r$-$1$}
\end{picture}&
$\af_i+\af_{i+1}+\cdots+\af_j$, $0\leq i\leq j<r$.\\\hline $B_r$&
\begin{picture}(76,15)
\multiput(5,5)(21,0){4}{\circle{2}}
\multiput(6,4)(0,2){2}{\line(1,0){19}}\put(48,5){\line(1,0){19}}
\put(30,5){$\ldots$}\put(13,3.5){\tiny$<$}
\put(4,-2){\tiny$0$}\put(25,-2){\tiny$1$}\put(42,-2){\tiny$r$-$2$}\put(63,-2){\tiny$r$-$1$}
\end{picture}&
$\af_i+\af_{i+1}+\cdots+\af_j$, $0\leq i\leq
j<r$,\\&&$2\af_0+\cdots+2\af_j+\af_{j+1}+\cdots+\af_k$, $0\leq j<k<r$.\\\hline
$C_r$&
\begin{picture}(76,15)
\multiput(5,5)(21,0){4}{\circle{2}}
\multiput(6,4)(0,2){2}{\line(1,0){19}}\put(48,5){\line(1,0){19}}
\put(30,5){$\ldots$}\put(13,3.5){\tiny$>$}
\put(4,-2){\tiny$0$}\put(25,-2){\tiny$1$}\put(42,-2){\tiny$r$-$2$}\put(63,-2){\tiny$r$-$1$}
\end{picture}&
$\af_i+\af_{i+1}+\cdots+\af_j$, $0\leq i\leq
j<r$,\\&&$\af_0+2\af_1+\cdots+2\af_j+\af_{j+1}+\cdots+\af_k$, $0\leq j\leq
k<r$.\\\hline

& & $\af_i+\af_{i+1}+\cdots+\af_j$, $0\leq i\leq j<r$,\\
$D_r$&
\begin{picture}(85,10)
\multiput(3,12)(0,-14){2}{\circle{2}} \put(4,-2){\line(2,1){12}}
\put(4,12){\line(2,-1){12}} \multiput(17,5)(21,0){4}{\circle{2}}
\put(18,5){\line(1,0){19}}\put(60,5){\line(1,0){19}} \put(42,5){$\ldots$}
\put(2,-9){\tiny$1$}\put(2,15){\tiny$0$}
\put(16,-2){\tiny$2$}\put(27,-2){\tiny$3$}\put(54,-2){\tiny$r$-$2$}\put(75,-2){\tiny$r$-$1$}
\end{picture}
&$\af_0+\af_2+\cdots+\af_j$, $2\leq j<r$,
\\&&$\af_0+\af_1+2\af_2+\cdots+2\af_j+\af_{j+1}+\cdots+\af_k$,
$2\leq j<k<r$.\\\hline
\end{tabular}
\end{center}
\end{table}

We have the following description of good Lyndon words. Calculations can be
found in Appendix \ref{A:Good Lyndon Words}.

\begin{prp} We have
\begin{enumerate}
\item The good Lyndon words for $\g$ of type $A_r$ are
\[
\{[i,\ldots,j]|0\leq i\leq j<r\}.
\]
\item The good Lyndon words for $\g$ of type $B_r$ are
\[
\{[i,\ldots,j]|0\leq i\leq j<r\}\cup\{[j,j-1,\ldots,0,0,\ldots,k-1,k]|0\leq
j<k<r\}.
\]
\item The good Lyndon words for $\g$ of type $C_r$ are
\[
\{[i,\ldots,j]|0\leq i\leq j<r\}\cup\{[j,\ldots,1,0,1,\ldots,k]|1\leq j<k\leq
r-1\}\cup\{[0,\ldots,j,1,\ldots,j]|1\leq j<r\}.
\]
\item The good Lyndon words for $\g$ of type $D_r$ are
\[
\{[0,2,\ldots,i]|2\leq i<r\}\cup\{[i,\ldots,j]|1\leq i\leq j\leq
r-1\}\cup\{[j,\ldots,1,0,2,\ldots,k]|1\leq j<k<r\}.
\]
\end{enumerate}
\end{prp}

We now list the root vectors associated to the good Lyndon words. Calculations
can be found in Appendix \ref{A:Root Vectors}

\begin{prp}
\begin{enumerate}
\item In type $A_r$,
\[
b^*[i,\ldots,j]=[i,\ldots,j],\;\;\;0\leq i\leq j<r.
\]

\item In type $B_r$:
\begin{align*}
b^*[i,\ldots,j]&=[i,\ldots,j],\;\;\;0\leq i\leq j<r\\
b^*[j,\ldots,0,0,\ldots,k]&=[2]_0[j,\ldots,0,0,\ldots,k],\;\;\;0\leq j<k<r.
\end{align*}

\item In type $C_r$:
\begin{align*}
b^*[i,\ldots,j]&=[i,\ldots,j],\;\;\;0\leq i\leq j<r,\\
b^*[j,\ldots,1,0,1,\ldots,k]&=[j,\ldots,1,0,1,\ldots,k],\;\;\;1\leq j<k<r,\\
b^*[0,\ldots,j,1,\ldots,j]&=q[0]\cdot([1,\ldots,j]*[1,\ldots,j]),\;\;\;1\leq j<r.
\end{align*}

\item In type $D_r$:
\begin{align*}
b^*[0]&=[0]\\
b^*[0,2,\ldots,i]&=[0,2,\ldots,i],\;\;\;2\leq i<r,\\
b^*[i,\ldots,j]&=[i,\ldots,j],\;\;\;1\leq i\leq j<r,\\
b^*_{[1, 0,2,\ldots,j]}&=[1,0,2,\ldots,j]+[0,1,2,\ldots,j],\;\;\;2\leq j<r,\\
b^*[j,\ldots,2,1,0,2,\ldots,k]&=[j,\ldots,2,1,0,2,\ldots,k]+[j,\ldots,2,0,1,2,\ldots,k],\;\;\;2\leq
j<k<r.
\end{align*}

\end{enumerate}
\end{prp}

\subsection{Good Lyndon Words in Type $E_8$}\label{SS:E8 Lyndon words}
Fix the following labeling of the nodes of the Dynkin diagram  for $E_8$.
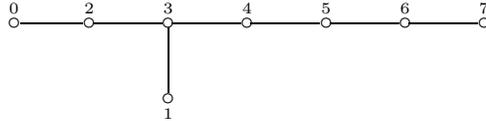
\[
\xy
(0,12)*{\mbox{\tiny0}};
(0,10)*{\circ};(10,10)*{\circ}**\dir{-};(10,12)*{\mbox{\tiny2}};
(10.5,10)*{};(20.5,10)*{\circ}**\dir{-};(20.5,-2)*{\mbox{\tiny1}};
(20.5,9.5)*{};(20.5,0)*{\circ}**\dir{-};(20.5,12)*{\mbox{\tiny3}};
(21,10)*{};(31,10)*{\circ}**\dir{-};(31,12)*{\mbox{\tiny4}};
(31.5,10)*{};(41.5,10)*{\circ}**\dir{-};(41.5,12)*{\mbox{\tiny5}};
(42,10)*{};(52,10)*{\circ}**\dir{-};(52,12)*{\mbox{\tiny6}};
(52.5,10)*{};(62.5,10)*{\circ}**\dir{-};(62.5,12)*{\mbox{\tiny7}};
\endxy
\]
We list here only the 12 good Lyndon words belonging to the set $\mathcal{E}$ in \cite[$\S8.3$]{klram2}:
\begin{align*}
&[6023145342302134567],\;[56023145345342302134567],\;[45623145342302134567],\\
&[3456023145342302134567],\;[13456023145342302134567],\;[23456023145342302134567],\\
&[323131456023145342302134567],\;[432131456023145342302134567],\\
&[543213456023145342342302134567],\;[6543213456023145342342302134567],\\
&[53423021345676451342302134567].
\end{align*}
The complete list of the 120 good Lyndon words for $E_8$ can be found in Appendix \ref{A:E8}.

\subsection{Good Lyndon Words and Root Vectors in Type $F_4$}\label{SS:TypeF}
We now calculate the Lyndon words and corresponding root vectors for $\g$ of type $F_4$. We choose the following labeling of the Dynkin diagram.
\[
\xy
(0,2)*{\circ};(10,2)*{\circ}**\dir{-};
(10.5,2)*{};(20.5,2)*{\circ}**\dir{=};(15.5,2)*{<};
(21,2)*{};(31,2)*{\circ}**\dir{-};
(0,0)*{\mbox{\tiny0}};(10,0)*{\mbox{\tiny1}};(20.5,0)*{\mbox{\tiny2}};(31,0)*{\mbox{\tiny3}};
\endxy
\]
Note that we have the opposite ordering as that in \cite{klram2}.

\begin{prp}\label{P:F4 Lyndon Words} The good Lyndon words for $F_4$ are given in the following table:
\begin{center}
\begin{tabular}{|c|c|}\hline
Height&Good Lyndon Words\\
\hline\hline
1&$[0], [1], [2],[3]$\\
2&$[01],[12],[23]$\\
3&$[012],[123],[112]$\\
4&$[0123],[1012],[1123]$\\
5&$[01012],[21123],[10123]$\\
6&$[010123],[210123]$\\
7&$[1210123],[2010123]$\\
8&$[12010123]$\\
9&$[112010123]$\\
10&$[2112010123]$\\
11&$[21012310123]$\\\hline
\end{tabular}
\end{center}
\end{prp}

\begin{prp} The root vectors for $F_4$ are given as follows:
\begin{equation*}
\begin{split}
b^*{[i,\ldots,j]}  =& [i,\ldots,j] \\
b^*{[112]}  =& [2]_0[112] \\
b^*{[1012]}  =& [1012] + [2]_0[0112] \\
b^*{[1123]}  =& [2]_0[1123] \\
b^*{[01012]}  =& [2]_0[01012] +[2]_0^2[00112] \\
b^*{[21123]}  =& [2]_0[21123] \\
b^*{[10123]}  =& [10123] + [2]_0[01123] \\
b^*{[010123]}  =& [2]_0[010123] + [2]_0^2[001123] \\
b^*{[210123]}  =& [210123] + [2]_0\left([201123]+[021123]\right) \\
b^*{[1210123]}  =& [1210123] + [2]_0\left([1021123]+[1201123]\right) \\
b^*{[2010123]}  =& [2]_0 \left([2010123]+ [0210123]\right) + [2]_0^2\left([2001123]+[0201123]+[0021123]\right) \\
b^*{[12010123]}  =& [2]_0\left([12010123]+[10210123]\right) + [01210123]\\ &+ [2]_0^2\left([12001123]+[10201123]+[10021123]\right)\\
               &  +  [2]_0([01201123]+[01021123])\\
b^*{[112010123]}  = & [2]_0[1]\cdot b^*[12020123] \\
b^*{[2112010123]}  = & [2]\cdot b^*[112010123]  \\
 b^*{[21012310123]}=& q[2]\cdot(b^*{[10123]}*b^*{[10123]})
\end{split}
\end{equation*}
\end{prp}

\begin{pff} For good Lyndon words of height at most 10, the result is obtained by direct calculation.
The calculation of the height 11 case is analogous to that for the long roots
in type $C_r$, see Proposition \ref{P:C root vectors} below.
\end{pff}

\subsection{The type $G_2$ case}\label{SS:TypeG} Fix the following labelling on the Dynkin diagram for $G_2$:
\[
\xy
{\ar@3{-}(0,2)*{};(10,2)*{}};(-.5,2)*\xycircle<2pt>{-};(10.5,2)*\xycircle<2pt>{-};(-.5,0)*{\mbox{\tiny0}};(10.5,0)*{\mbox{\tiny1}};(5,2)*{<};
\endxy
\]
Below we list the good Lyndon words, and associated root vectors:

\begin{prp} The good Lyndon words for $\Gm$ of type $G_2$ are
\[
[0],[1],[01],[001],[0001],[00101].
\]
\end{prp}

\begin{prp} The root vectors for $G_2$ are as follows:
\begin{align*}
b^*[0]&=[0],\\b^*[1]&=[1],\\b^*[01]&=[01],\\b^*[001]&=[2]_0[001]\\b^*[0001]&=[2]_0[3]_0[0001],\\b^*[00101]&=[2]_0[3]_0[00101]+[2]_0[3]_0[2]_1[00011].
\end{align*}
\end{prp}

We note here that the Lyndon words and associated root vectors agree whether we
read from right-to-left or from left-to-right, cf. \cite[$\S5.5.4$]{lec}.

\section{Construction of the Cuspidal Representations} \label{A:SegReps}

Fix the (right) Lyndon ordering on $\G$ as in Section \ref{S:IdentGoodLyndon}.
Recall that $\ui_l$ denotes the reading of a good Lyndon word $l$. That is
$l=[\ui_l]$. Throughout this section, we will need the converse to \cite[Lemma 6.4]{klram2}. The proof is very similar to \cite[Lemma 6.6]{klram2}.

\begin{lem}\label{L:IrreduciblilityCriterion}
Let $V\in H(\Gm;\bt)$, and assume that $\ch V=b^*_l$ for some $l=l(\bt)\in\GL$. Then, $V$ is irreducible.
\end{lem}

\begin{pff}
By Theorem \ref{T:minbg}, all composition factors of $V$ have lowest weight $g\in\G$ for $g>l$. On the other hand, all composition factors have lowest weight belonging to $\G_\bt$ so, by Corollary \ref{C:MaxGoodWord}, $V=L^{\oplus k}$ for some simple module $L$. The result now follows because $\{ b^*_g\,|\,g\in\G\}$ is an $\A$ basis of $\U^*_\A$.
\end{pff}

\subsection{Type $A_r$}\label{SS:A}
Let $l=[i,\ldots,j]$, $0\leq i\leq j<r$. We have $b^*_l = [i \dubc j] $. Define
$\va_l= \bbF . v_{0}$ where $v_0$ has degree $0$.  Set $e(\ui)\va_l =
\delta_{\ui,\ui_l}\va_l$, $\phi_s v_0 = 0$ and $y_s v_0 = 0$ for all admissible
$s$.  This is the trivial representation of $H(\Gm;\nu)$ and clearly satisfies
\eqref{E:Relations1}-\eqref{E:Braid Relations} and $\ch\va_l=b^*_l$.

\subsection{Type $B_r$}\label{SS:B}
The case $l=[i,\ldots,j]$, $0\leq i\leq j<r$ is the trivial representation as
in type $A_r$.

Let $l=[j,\ldots,0,0,\ldots,k]$, $0\leq j<k<r$. Then,
$b^*_l=(q+q^{-1})[j,\ldots,0,0,\ldots,k]$. Set $\va_l=\bbF v_1\oplus\bbF
v_{-1}$, where $\deg v_i=i$ for $i=\pm1$. Define
$e(\ui)\va_1=\dt_{\ui,\ui_l}\va_l$. Set $y_rv_1=0$ for all $s$, for $s\neq
j+1$, set $\phi_s v_1=0$, and define $\phi_{j+1}v_1=v_{-1}$. Set
$\phi_rv_{-1}=0$ for all $s$, for $s\neq j+1,j+2$, set $\phi_rv_{-1}=0$, and
set $y_{j+1}v_{-1}=-v_1$ and $y_{j+2}v_{-1}=v_1$. We leave it as an easy
exercise to the reader to check that this satisfies
\eqref{E:Relations1}-\eqref{E:Braid Relations} and $\ch\va_l=b^*_l$

\subsection{Type $C_r$}\label{SS:C}
For $l\neq [0,\ldots,j,1,\ldots,j]$, $\va_l$ is the trivial representations and
may be computed as in type $A_r$.

Assume $l=[0,\ldots,j,1,\ldots,j]$. Then
$b^*_l=q[0]([1,\ldots,j]*[1,\ldots,j])$. Let $\bt=\af_1+\cdots+\af_j$, and
consider the $H(\Gm;\af_0,2\bt)$ module
$\va_{\af_0}\boxtimes(\ind_{\bt,\bt}^{2\bt}\va_{\bt}\boxtimes\va_\bt)\{1\}$.
Extend this to a $H(\Gm;2\bt+\af_0)$ module by insisting that $\phi_1$ acts as
0, and $e(\ui)$ acts as 0 if $i_1\neq0$. It is very easy to check that this is
the desired cuspidal representation, $\va_l$, cf. \cite[$\S8.6$]{klram2}.

\subsection{Type $D_r$}\label{SS:D}
For $l\neq [1,0,2,\ldots,k]$, $1\leq j<k<r$, $\va_l$ is the trivial
representation and can be computed as in type $A_r$.

Assume $l=[j,\ldots,1,0,2,\ldots,k]$. Define $\va_l=\bbF v_0\oplus\bbF w_0$,
where $v_0$ and $w_0$ have degree $0$. Define
\[
e(\ui)\va_l=\dt_{\ui,\ui_l}\bbF v_0 +\dt_{\ui, s_j\cdot\ui_l}\bbF w_0.
\]
Define $y_r\va_l=0$. For $r\neq j$, define $\phi_r\va_l=0$ and set $\phi_j
v_0=w_0$. It is elementary to check that this is indeed a representation and
$\ch\va_l=b^*_l$.

\subsection{Type $E_8$}\label{SS:E}
We simply note here that in our ordering all Lyndon words for type $E_8$ are
homogeneous in the sense of \cite{klram} and the corresponding cuspidal
representations can be computed using \cite[Theorems 3.6,3.10]{klram}. The 12
outstanding cases from \cite{klram2} are listed in subsection \ref{SS:E8 Lyndon
words} and are evidently homogeneous. An entire list of the good Lyndon words
for $E_8$ can be found in Appendix \ref{A:E8}.

\subsection{Type $F_4$}\label{SS:F}
We choose the following partial ordering on $I$: $0\to1\to2\to3$, see
\eqref{E:finite type Q polynomials}.

\begin{enumerate}
\item $l=[i,\ldots,j]$, $0\leq i\leq j\leq 3$.

Constructed exactly as in the type $A$ case.

\item $l=[112],[1123],\mbox{ or }[21123].$

Constructed exactly as in the type $B$ case.

\item $l=[1012],\mbox{ or }[01012]$.

These are constructed exaclty as in the type $C$ case in \cite{klram2}. Indeed, we have
\[
b^*[1012]=([1]*[01])[2],\andeqn b^*[01012]=q([01]*[01])[2].
\]
For example, let $\bt=\af_0+\af_1$ and define the $H(\Gm;2\bt,\af_2)$-module
$$V=(\ind_{\bt,\bt}^{2\bt}\va_{\bt}\boxtimes\va_{\bt})\boxtimes\va_{\af_2}\{1\}.$$
Extend the action to $H(\Gm;2\bt+\af_2)$ by insisting that $\phi_4$ acts as 0
and $e(\ui)$ acts as 0 if $i_5\neq2$. As in \cite{klram2}, the only relation
that is not obvious is \ref{E:Quadratic Relations}, which follows since
$y_4^2-y_5$ acts as 0 on the module above. Then, $\va_{[01012]}=V$ is the
desired cuspidal representation.

\item $l=[10123]$ or $[010123]$.

In either case, let $\bt=|l|$. Define the $H(\Gm;\bt,\af_3)$-module
$V=\va_\bt\boxtimes\va_{\af_3}$. As above, we may extend this to a
$H(\Gm;\bt+\af_3)$-module by insisting that $\phi_r$ acts as 0 and $e(\ui)$
acts as 0 if $i_{r+1}\neq3$, where $r=4,\mbox{ or }5$ as appropriate. To check
relation \ref{E:Quadratic Relations} it is enough to observe that $y_r-y_{r+1}$
acts as 0 on $V$ (actually, each both $y_r$ and $y_{r+1}$ act as 0). Hence,
$\va_\bt=V$ is the desired cuspidal representation.

\item $l=[210123]$.

Let $\bt=|l|-\af_2$. Define a graded vector space $V=W\oplus U$, where
$W\cong\va_{\af_2}\boxtimes\va_{\bt}$ as a $H(\Gm;\af_2,\bt)$-module, and
$U=U[1]\oplus U[-1]$ is 2-dimensional with basis $\{u_1,u_{-1}\}$. Fix an
\emph{weight} basis $\{w_0,w_1,w_{-1}\}$ for $W$. That is,
$e(\ui)w_0=\dt_{\ui,\ui_l}w_0$ and $w_0$ has degree 0, $w_1=\phi_2w_0$, and
$w_{-1}=\phi_3w_1$. It follows from \eqref{E:Braid Relations} that
$\phi_1w_{-1}=w_0$.

The following defines an action of $H(\Gm;\af_2+\bt)$ on $V$:
\begin{itemize}
\item $\phi_1w_0=0$, $\phi_1w_1=u_1$, $\phi_1w_{-1}=u_{-1}$;
\item $\phi_2u_{-1}=0$;
\item $e(\ui)$ acts as 0 on $W$ if $i_1\neq2$.
\end{itemize}
Indeed, from \eqref{E:Relations4} we are forced to define
$$e(\ui)u_{\pm1}=e(\ui)\phi_1w_{\pm1}=\phi_1e(s_1\cdot\ui)w_{\pm1}=\begin{cases}u_{\pm1}&\mbox{if }[\ui]=[021123],\\0&\mbox{otherwise.}\end{cases}$$
Using \eqref{E:Relations3}-\eqref{E:Relations9}, we must set $y_ru_1=0$ for
$1\leq r\leq 6$. For example,
\[
y_1u_1=y_1\phi_1\phi_2w_0=\phi_1y_2\phi_2w_0=\phi_1\phi_2y_3w_0=0.
\]
Also, we define $y_ru_{-1}=0$ if $r\neq 3,4$, and
\[
y_3u_{-1}=y_3\phi_1\phi_3\phi_2w_0=\phi_1(\phi_3y_4-1)\phi_2w_0=-u_1.
\]
Similarly, $y_4u_{-1}=u_1$.

Using \eqref{E:Relations7}, we define $\phi_3u_1=\phi_1w_{-1}=u_{-1}$, and
$\phi_4u_{\pm1}=\phi_5u_{\pm1}=0$. The relation \eqref{E:Quadratic Relations}
forces $\phi_1u_1=w_1$, $\phi_1u_{-1}=w_{-1}$, and $\phi_3u_{-1}=0$. Using
\eqref{E:Braid Relations} we define
$$\phi_2u_1=\phi_2\phi_1\phi_2w_0=\phi_1\phi_2\phi_1w_0=0.$$

We need to show that the actions of $\phi_1$ and $\phi_2$ are consistent with
relations \eqref{E:Relations1}-\eqref{E:Braid Relations}. As explained above,
relations \eqref{E:Relations1}-\eqref{E:Relations7} follow by definition, as do
the relations \eqref{E:Relations8}-\eqref{E:Quadratic Relations} for the action
of $\phi_1$.

We will postpone checking \eqref{E:Braid Relations} until we have checked the
action of $\phi_2$ on $U$. To check relations \eqref{E:Relations8} and
\eqref{E:Relations9} we need only consider the nontrivial cases $r=3,4$.
Indeed, we compute
\begin{align*}
y_4\phi_2u_{-1}&=y_4\phi_2(\phi_1\phi_3\phi_2)w_0&\\
    &=\phi_2\phi_1(\phi_3y_3+1)\phi_2w_0&\mbox{by \eqref{E:Relations8} in $W$}\\
    &=\phi_2\phi_1\phi_3\phi_2y_2w_0+\phi_2\phi_1\phi_2w_0&\\
    &=\phi_1\phi_2\phi_1w_0&\mbox{by \eqref{E:Braid Relations} in $W$}\\
    &=0,&
\end{align*}
since $\phi_1w_0=y_2w_0=0$. A similar computation with $r=3$ gives
\eqref{E:Relations9}. To check relation \eqref{E:Quadratic Relations} we need
only observe $y_1u_{-1}=y_2u_{-1}=0$. Finally, the last nontrivial relation is
\begin{align*}
\phi_1\phi_2u_{-1}&=\phi_1\phi_2\phi_1\phi_3\phi_2w_0
    =\phi_2\phi_1\phi_2\phi_3\phi_2w_0\\
    &=\phi_2\phi_1(\phi_3\phi_2\phi_3-1)w_0
    =0.
\end{align*}
One has $\ch V=b^*[210123]$. Hence, $\va_l=V$ is the desired representation.

\item $l=[1210123]$.

Let $\bt=|l|-\af_1$ and define the $H(\Gm;\af_1,\bt)$-module
$V=\va_{\af_1}\boxtimes\va_\bt$. Extend this to an action of $H(\Gm;\bt+\af_1)$
by insisting that $\phi_1$ acts as 0 and $e(\ui)$ acts as 0 if $i_1\neq 1$.
Again, the only thing nontrivial to check is \eqref{E:Quadratic Relations}
which follows since $y_1^2-y_2$ acts as 0 on $V$ (actually, both $y_1$ and
$y_2$ act as 0). Then, $\va_l=V$ is the desired representation.

\item $l=[2010123]$.

Let $\bt=|l|-\af_2$. Define the graded vector space $V=W\oplus U$, where
$W\cong \va_{\af_2}\boxtimes\va_\bt$ as a $H(\Gm;\af_2,\bt)$-module and
$U=U[2]\oplus U[1]\oplus U[0]\oplus U[-1]\oplus U[-2]$ is 10-dimensional with
basis $\{u_2^{1},u_2^2, u_1,u_0^1, u_0^2,u_{-0}^{1}, u_{-0}^2, u_{-1},
u_{-2}^1, u_{-2}^2\}$. Fix a weight basis
$\{w_2,w_1,w_0,w_{-0},w_{-1},w_{-2}\}$ for $W$. That is,
$e(\ui)w_1=\dt_{\ui,\ui_l}w_1$, $\deg w_1=1$, $w_2=\phi_3w_1$,
$w_{0}=\phi_2w_2$, $w_{-0}=\phi_4w_2$, $w_{-2}=\phi_4w_0=\phi_2w_{-0}$, and
$w_{-1}=\phi_3w_{-2}$.

The following defines an action of $H(\Gm,\bt+\af_2)$ on $V$:
\begin{itemize}
\item $\phi_1w_{\pm1}=u_{\pm1}$;
\item For $i\in\{2,0,-0,-2\}$, $\phi_1w_i=u_i^1$;
\item $\phi_2u_1=\phi_2u_{-1}=0$;
\item For $i\in\{2,0,-0,-2\}$, $\phi_2u_i^1=u_i^2$;
\item $\phi_3u_{-0}^2=\phi_3u_{-2}^2=0$
\item $e(\ui)$ acts as 0 on $W$ if $i_1\neq2$.
\end{itemize}
The remaining relations are now forced. By \eqref{E:Relations4} we have
\begin{itemize}
\item $e(\ui)u_{\pm1}=\begin{cases}u_{\pm1}&\mbox{if }\ui=[0210123],\\0&\mbox{otherwise;}\end{cases}$
\item For $i\in\{2,0,-0,-2\}$, $e(\ui)u_i^1=\begin{cases}u_i^1&\mbox{if }\ui=[0201123],\\0&\mbox{otherwise;}\end{cases}$
\item For $i\in\{2,0,-0,-2\}$, $e(\ui)u_i^2=\begin{cases}u_i^2&\mbox{if }\ui=[0021123],\\0&\mbox{otherwise.}\end{cases}$
\end{itemize}
We now use \eqref{E:Relations6} and \eqref{E:Relations8}-\eqref{E:Relations9}
to define the action of $y_1,\ldots,y_7$ on $U$. Since $y_4,\ldots,y_7$ commute
with $\phi_1$ and $\phi_2$, their actions are determined by $W$. As an example,
we compute the action of $y_1$ on $U$ below. The action of $y_2$ and $y_3$ can
be worked out similarly.
\[
\begin{array}{lll}y_1u_2^1=\phi_1y_2w_2,&y_1u_2^2=\phi_2(\phi_1y_2-1)w_2,&y_1u_1=\phi_1y_2w_1,\\
y_1u_0^1=\phi_1y_2w_0,&y_1u_0^2=\phi_2(\phi_1y_2-1)w_0,&y_1u_{-0}^1=\phi_1y_2w_{-0},\\
y_1u_{-0}^2=\phi_2(\phi_1y_2-1)w_{-0},&y_1u_{-1}=\phi_1y_2w_{-1},&y_1u_{-2}^1=\phi_1y_2w_{-2},\\
y_1u_{-2}^2=\phi_2(\phi_1y_2-1)w_{-2}.&&\end{array}
\]

Next, to define the action of $\phi_1,\ldots,\phi_6$ on $U$, we note that since
$\phi_3,\ldots,\phi_6$ commute with $\phi_1$, their actions on $u_2^1$, $u_1$,
$u_0^1$, $u_{-0}^1$, $u_{-1}^1$ and $u_{-2}^1$ are determined by $W$.
Additionally, since $\phi_4,\phi_5,\phi_6$ commute with $\phi_1$ and $\phi_2$,
their action on $u_2^2$, $u_0^2$, $u_{-0}^2$, and $u_{-2}^2$ are determined by
$W$. The remaining calculations are given below and can be worked out by
rewriting the $u$'s in the form $\phi_\sm w_1$.
\[
\begin{array}{llll}\phi_1u_1=w_1,&\phi_1u_{-1}=w_{-1},&&\\
\phi_1u_2^1=w_2,&\phi_1u_0^1=w_0,&\phi_1u_{-0}^1=w_{-0},&\phi_1u_{-2}^1=w_{-2}\\
\phi_1u_2^2=u_1^2,&\phi_1u_0^2=0,&\phi_1u_{-0}^2=u_0^2,&\phi_1u_{-2}^2=0\\
\phi_2u_2^2=u_2^1,&\phi_2u_0^2=u_0^1,&\phi_2u_{-0}^2=u_{-0}^1,&\phi_2u_{-2}^2=u_{-2}^1\\
\phi_3u_2^2=0,&\phi_3u_0^2=0,&&
\end{array}
\]

We now have to check that the actions of $\phi_2$ on $u_{\pm1}$ and $\phi_3$ on
$u_{-0}^2,u_{-2}^2$ are consistent with the relations. Indeed, in the case
$\phi_2u_1=0$, the only nontrivial relations to check are \eqref{E:Quadratic
Relations} and \eqref{E:Braid Relations}. We have for \eqref{E:Quadratic
Relations},
\begin{align*}
\phi_2^2u_1&=Q_{21}(y_2,y_3)u_1\\
    &=(y_3^2-y_2)\phi_1w_1\\
    &=\phi_1(y_3^2-y_1)w_1\\
    &=0.
\end{align*}
For the braid relations, we have
$$\phi_1\phi_2u_1=\phi_1\phi_2\phi_1w_1=\phi_2\phi_1\phi_2w_1=0,$$ and
$$\phi_2\phi_3\phi_2u_1=\phi_3\phi_2\phi_3\phi_1w_1=\phi_3\phi_2\phi_1\phi_3w_1=\phi_3u_2^2=0.$$
We now check that $\phi_2u_{-1}=0$ is consistent with the relations. Indeed,
one verifies that
\begin{align*}
\phi_2^2u_{-1}&=Q_{21}(y_2,y_3)u_{-1}=(y_3^2-y_2)u_{-1}=0
\end{align*}
For the braid relations, we have
$$\phi_1\phi_2u_{-1}=\phi_1\phi_2\phi_1w_{-1}=\phi_2\phi_1\phi_2w_{-1}=0,$$
and
\begin{align*}
\phi_2\phi_3\phi_2u_{-1}&=\phi_3\phi_2\phi_3u_{-1}\\
    &=(\phi_3\phi_2\phi_1)\phi_3^2w_{-2}\\
    &=(\phi_3\phi_2\phi_1)Q_{01}(y_3,y_4)w_{-2}\\
    &=(\phi_3\phi_2\phi_1)(y_3-y_4)\phi_4\phi_2\phi_3w_1\\
    &=\phi_3\phi_2\phi_1(\phi_4\phi_3+\phi_2\phi_3)w_1\\
    &=\phi_3(u_{-0}^2+u_0^2)=0.
\end{align*}

We now check the action of $\phi_3$. Indeed, for \eqref{E:Quadratic Relations},
a calculation gives
$$\phi_3^2u_{-0}^2=Q_{21}(y_3,y_4)u_{-0}^2=(y_4^2-y_3)u_{-0}^2=0.$$
Similarly, $\phi_3^2u_{-2}^2=0$. For \eqref{E:Braid Relations}, we need only
calculate
\begin{align*}
\phi_2\phi_3u_{-0}^2&=\phi_2\phi_3\phi_2u_{-0}^1\\
    &=\phi_3\phi_2\phi_3(\phi_1\phi_4\phi_3)w_1\\
    &=\phi_3\phi_2\phi_1(\phi_4\phi_3\phi_4-1)w_1\\
    &=-\phi_3\phi_2\phi_1w_1\\
    &=-\phi_3\phi_2u_1=0.
\end{align*}
Similarly, we have $\phi_2\phi_3u_{-2}^1=0$.

We have $\ch V=b^*[2010123]$. Hence, $\va_l=V$ is the desired representation.

\item $l=[12010123]$.

Let $\bt=|l|-\af_1$. Define the graded vector space $V=(W\oplus U)\oplus Z$,
where $W\oplus U\cong \va_{\af_1}\boxtimes\va_\bt$ as a
$H(\Gm;\af_1,\bt)$-module and has a basis as described in the previous case and
$Z=Z[1]\oplus Z[0]\oplus Z[-1]$ is 5-dimensional with basis
$\{z_1^1,z_1^2,z_0,z_{-1}^1,z_{-1}^2\}$.

The following defines an action of $H(\Gm;\bt+\af_1)$ on $V$:
\begin{itemize}
\item $\phi_1u_0^1=z_1^1$;
\item $\phi_1u_0^2=z_1^2$;
\item $\phi_1u_{-1}=z_0$;
\item $\phi_1u_{-2}^1=z_{-1}^1$;
\item $\phi_1u_{-2}^2=z_{-1}^2$;
\item $\phi_1$ acts as 0 on the remaining basis vectors of $W\oplus U$;
\item $e(\ui)$ act as 0 on $W\oplus U$ if $i_1\neq 1$.
\end{itemize}

We now determine the remaining actions of $H(\Gm;\bt+\af_1)$ on $Z$. Indeed,
note that $y_3,\ldots,y_8$ commute with $\phi_1$, so their actions are
determined by $W\oplus U$. To calculate the action of $y_1$ and $y_2$, note
that as operators on $Z$, $y_1\phi_1=\phi_1y_2$ and $y_2\phi_1=\phi_1y_1$ so
their action is determined by $U$. In particular, $y_2$ acts as 0 on $Z$ since
$y_1$ acts as 0 on $U$. Additionally, a calculation gives
\[
\begin{matrix}y_1z_1^1=-\phi_1u_2^1=0,&y_1z_1^2=\phi_1u_2^2=0,&y_1z_0=\phi_1u_1=0,\\
y_1z_{-1}^1=-\phi_1u_{-0}^1=0,&y_1z_{-1}^2=-\phi_1u_{-0}^2=0.&\end{matrix}
\]
Next observe that the action of $\phi_3,\ldots,\phi_8$ on $Z$ are determined by
$W\oplus U$. We calculate
\begin{align}\label{E:PhionZ}
\begin{matrix}\phi_1z_1^1=-u_2^1,&\phi_1z_1^2=-u_2^2,&\phi_1z_0=0,&\phi_1z_{-1}^1=-u_{-0}^1,&\phi_1z_{-1}^2=-u_{-0}^2\end{matrix}
\end{align}
and $\phi_2$ acts as 0 on $Z$.

It remains to check the consistency of this action with the relations. The only
relations which are not obvious are \eqref{E:Quadratic Relations} and
\eqref{E:Braid Relations} for $\phi_1$.

To check \eqref{E:Quadratic Relations} on $W$ it is enough to check that
$\phi_1^2w_1=Q_{12}(y_1,y_2)w_1=0$ which is obvious. Many of the quadratic
relations for the action of $\phi_1$ on $U$ are contained in \eqref{E:PhionZ}
above. The remaining calculation are outlined below.
\[
\begin{matrix}\phi_1^2u_2^1=Q_{10}(y_1,y_2)u_2^1=0,&\phi_1^2u_2^2=\phi_3\phi_1^2u_2^1=0,&\phi_1^2u_1=Q_{10}(y_1,y_2)u_1=0,\\
    \phi_1^2u_{-0}^1=Q_{10}(y_1,y_2)u_{-0}^1=0&\phi_1^2u_{-0}^2=\phi_3\phi_1^2u_{-0}^1=0.&\end{matrix}
\]
Relation \eqref{E:Quadratic Relations} for the action of $\phi_1$ on $Z$ is now
obvious.

To check \eqref{E:Braid Relations} we need to show that $\phi_1\phi_2z=0$ for
all $z\in Z$. This calculation, however, is trivial. For example,
$$\phi_1\phi_2z_1^1=\phi_1\phi_2\phi_1u_{0}^1=\phi_2\phi_1\phi_2u_0^1=\phi_2\phi_1\phi_2^2w_0=\phi_2\phi_1w_0=0.$$
We have $\ch V=b^*[12010123]$. Hence, $\va_l=V$ is the desired cuspidal representation.

\item $l=[112010123]$.

Let $\bt=|l|-\af_1$. Define the graded vector space $V=W\{1\}\oplus W\{-1\}$,
where $W\cong\va_{\af_1}\boxtimes\va_\bt$ as a $H(\Gm;\af_1,\bt)$-module. For
each $w\in W$  write $w\{\pm1\}\in W\{\pm1\}$ for the corresponding vector.

The following defines an action of $H(\Gm;\bt+\af_1)$ on $V$:
\begin{itemize}
\item $\phi_1w\{1\}=w\{-1\}$ for $w=e(1,1,\ldots)w$ and $\phi_1 w\{1\}=0$ if $w=e(1,0,\ldots)w$;
\item $\phi_1$ acts as 0 on $W\{-1\}$;
\item $e(\ui)$ acts as 0 on $V$ unless $i_1=1$.
\end{itemize}

Once again, it is straightforward to see that this is an
$H(\Gm;\bt+\af_1)$-module. Indeed, the only relation to check is
\eqref{E:Quadratic Relations}. For $v\in V$, note that either $v=e(11\ldots)v$
or $v=e(10\ldots)v$. Hence
\[
\phi_1^2v=\begin{cases}0&\mbox{if }v=e(11\ldots)v;\\
(y_2-y_1)v&\mbox{if }v=e(10\ldots)v.\end{cases}
\]
The result now follows since both $y_1$ and $y_2$ act as 0 on $V$. We have $\ch
V=b^*[112010123]$, so $\va_l=V$ is the desired representation.

\item $l=[2112010123]$.

Let $\bt=|l|-\af_2$ and define $V=\va_{\af_2}\boxtimes\va_\bt$ as a
$H(\Gm;\af_2,\bt)$-module. Extend this to an $H(\Gm;\bt+\af_2)$-module by
insisting that $\phi_1$ acts as 0 and $e(\ui)$ acts as 0 unless $i_1=2$. The
only relation that is nontrivial to check is \eqref{E:Quadratic Relations},
which follows since $y_2^2-y_1$ acts as 0 on $V$. Hence $\va_l=V$ is the
desired representation.

\item $l=[21012310123]$.

Let $\bt=\af_0+2\af_1+\af_2+\af_3$. Consider the $H(\Gm;\af_2,2\bt)$ module
$$V=\va_{\af_2}\boxtimes(\ind_{\bt,\bt}^{2\bt}\va_{\bt}\boxtimes\va_\bt)\{1\}.$$
Extend this to a $H(\Gm;2\bt+\af_2)$ module by insisting that $\phi_1$ acts as
0 and $e(\ui)$ acts as 0 if $i_1\neq 2$. As in the case of the long roots of
type $C$, clearly the relations for $H(\Gm;\af_2,2\bt)$ are satisfied. The only
new relation which is not obvious is \eqref{E:Quadratic Relations}, which
follows since $y_2^2-y_1$ acts as 0 on $V$. Hence we have constructed a module
$\va_l=V$ with character $b^*[21012310123]$.
\end{enumerate}

\appendix

\section{Calculations}\label{S:Appendix}

\subsection{Good Lyndon Words}\label{A:Good Lyndon Words}

\begin{prp} We have
\begin{enumerate}
\item The Good Lyndon words for $\g$ of type $A_r$ are
\[
\{[i,\ldots,j]|0\leq i\leq j<r\}.
\]
\item The good Lyndon words for $\g$ of type $B_r$ are
\[
\{[i,\ldots,j]|0\leq i\leq j<r\}\cup\{[j,j-1,\ldots,0,0,\ldots,k-1,k]|0\leq
j<k<r\}.
\]
\item The good Lyndon words for $\g$ of type $C_r$ are
\[
\{[i,\ldots,j]|0\leq i\leq j<r\}\cup\{[j,\ldots,1,0,1,\ldots,k]|1\leq j<k\leq
r-1\}\cup\{[0,\ldots,j,1,\ldots,j]|1\leq j<r\}.
\]
\item The good Lyndon words for $\g$ of type $D_r$ are
\[
\{[0,2,\ldots,i]|2\leq i<r\}\cup\{[i,\ldots,j]|1\leq i\leq j\leq
r-1\}\cup\{[j,\ldots,1,0,2,\ldots,k]|1\leq j<k<r\}.
\]
\end{enumerate}
\end{prp}

\begin{pff}
Proceed by induction on the $\height(\bt)$. In all types, $\af_i\in\Pi$ implies
$l(\af_i)=[i]$.

\begin{enumerate}
\item For $\beta = \alpha_{i} + \cdots + \alpha_{j}$, we have
$$C\left(\beta\right) = \left\{\left(\alpha_{i} + \cdots + \alpha_{k} , \alpha_{k+1} + \cdots + \alpha_{j}\right) : j>k \geq i\right\}.$$
By induction, we assume
$$l\left(\alpha_{k+1} + \cdots + \alpha_{j}\right) = \left[k+1 \dubc j\right]\andeqn l\left(\alpha_{i} + \cdots + \alpha_{k}\right)= \left[i \dubc k\right].$$
Thus, $l\left(\beta\right) = \min\left\{\left[i \dubc k , k+1 \dubc j\right] : j>k \geq i\right\}  = \left[i \dubc j\right]$ completing our induction.

\item For $\beta = \alpha_{i} + \cdots \alpha_{j}$, we repeat the argument for type A to obtain $l\left(\beta\right) = [i , \ldots , j]$.

We now calculate $l(\bt)$ for $\beta = 2\alpha_{0} + \cdots + 2\alpha_{j} + \alpha_{j+1} + \cdots + \alpha_{k}$. We have
\begin{align*}
C\left(\beta\right) =& \left\{\left(2\alpha_{0} + \cdots + 2\alpha_{j} + \alpha_{j+1} + \cdots + \alpha_{i} , \alpha_{i+1} + \cdots + \alpha_{k}\right) | k>i>j\right\}\\ &\cup \left\{\left(2\alpha_{0} + \cdots + 2\alpha_{i} + \alpha_{i+1} + \cdots + \alpha_{j} , \alpha_{i+1} + \cdots + \alpha_{k}\right) | k>j>i\right\}\\
&\cup \left\{\left(\alpha_{i+1} + \cdots + \alpha_{j} , 2\alpha_{0} + \cdots + 2\alpha_{i} + \alpha_{i+1} + \cdots + \alpha_{k}\right) | k>j>i\right\}\\ &\cup \left\{\left(\alpha_{0} + \cdots + \alpha_{j} , \alpha_{0} + \cdots + \alpha_{k}\right)\right\}.
\end{align*}
Recall that
$l\left(\alpha_{i+1} + \cdots + \alpha_{k}\right) = [i+1 , \ldots , k]$, $l\left(\alpha_{i+1} + \cdots + \alpha_{j}\right) = [i+1 , \ldots , j]$, $l\left(\alpha_{0} + \cdots + \alpha_{j}\right) = [0 , \ldots , j]$ and $l\left(\alpha_{0} + \cdots + \alpha_{k}\right) = [0 , \ldots , k]$.

Our base case is $2\alpha_{0} + \alpha_{1} = \beta$.
Here, the first three sets which constitute $C(\bt)$ are empty and $l\left(\beta\right) = l\left(\alpha_{0}\right)l\left(\alpha_{0} + \alpha_{1}\right) = [0 , 0 , 1]$.

Assume by induction on the height of $\bt$ that
$$l\left(2\alpha_{0} + \cdots + 2\alpha_{j} + \alpha_{j+1} + \cdots + \alpha_{i}\right) = [j , \ldots , 0 , 0 , \ldots , i ],\;\;\;j<i<k,$$
$$l\left(2\alpha_{0} + \cdots + 2\alpha_{i} + \alpha_{i+1} + \cdots + \alpha_{j}\right) = [i , \ldots , 0 , 0 , \ldots , j ],\;\;\;i<j,$$
$$l\left(2\alpha_{0} + \cdots + 2\alpha_{i} + \alpha_{i+1} + \cdots + \alpha_{k}\right) = [i , \ldots , 0 , 0 , \ldots , k ],\;\;\;i<j.$$
Then,
\begin{align*}
l\left(\beta\right) =& \min\left\{[j \dubc 0 , 0 \dubc i , i+1 \dubc k] | i>j \right\}\\ &\cup \left\{[i \dubc 0 , 0 \dubc j , i+1 \dubc k] | j>i \right\} \\
&\cup \left\{[i+1 \dubc j , i \dubc 0 , 0 \dubc k] | j>i \right\}\\ &\cup \left\{[0 \dubc j , 0 \dubc k]\right\} \\
=& [j , \ldots , 0 , 0 , \dubc k]
\end{align*}
completing our induction.

\item For $\beta = \alpha_{i} + \cdots + \alpha_{j}$, we repeat the argument for type A to obtain $l\left(\beta\right) = [i , \ldots , j]$.

The next cases are somewhat more subtle. Observe for the base case that $$C\left(\alpha_{0} + 2\alpha_{1}\right) = \left\{\left(\alpha_{0} + \alpha_{1} , \alpha_{1}\right)\right\}$$ so that $l\left(\alpha_{0} + 2\alpha_{1}\right) = [0 , 1 , 1]$. Also,
$$C(\af_0+2\af_1+\af_2)=\{(\af_0+2\af_1,\af_2),(\af_0+\af_1,\af_1+\af_2),(\af_1,\af_0+\af_1+\af_2)\}.$$
Evidently, this gives $l(\af_0+2\af_1+\af_2)=[1,0,1,2]$.

Assume that $\bt=\af_0+2\af_1+\cdots+2\af_j$, and we have shown that for $i<k\leq j$,
$$l(\af_0+2\af_1+\cdots+2\af_i)=[0,\ldots,i,1,\ldots,i]$$ and $$l(\af_0+2\af_1+\cdots+2\af_i+\af_{i+1}+\cdots+\af_k)=[i,\ldots,1,0,1,\ldots,k].$$
Observe
\begin{align*}
C(\bt)=&\{(\af_0+2\af_1+\cdots+2\af_i+\af_{i+1}+\cdots+\af_j,\af_{i+1}+\cdots+\af_j)|1\leq i<j\}\\&\cup\{(\af_0+\cdots+\af_j,\af_1,\ldots,\af_j)\}.
\end{align*}
Thus,
\begin{align*}
l(\bt)=&\min\{[i,\ldots,1,0,1,\ldots,j,i+1,\ldots,j]|1\leq i<j\}\\&\cup\{[0,\ldots,j,1,\ldots,j]\}\\=&[0,\ldots,j,1,\ldots,j].
\end{align*}

Finally, assume $k>j$ and $\bt=\af_0+2\af_1+\cdots+2\af_j+\af_{j+1}+\cdots+\af_k$. Assume further that for all $j<i<k$
$$l(\af_0+2\af_1+\cdots+2\af_j+\af_{j+1}+\cdots+\af_i)=[j,\ldots,1,0,1,\ldots,i],$$
and assume that for either $i<j$ and $m\leq k$, or $i=j$ and $m<k$ that
$$l(\af_0+2\af_1+\cdots+2\af_i+\af_{i+1}+\cdots+\af_m)=[i,\ldots,1,0,1,\ldots,m].$$
We have
\begin{align*}
C\left(\beta\right) =& \left\{\left(\alpha_{0} + 2\alpha_{1} + \cdots + 2\alpha_{j} , \alpha_{j+1} + \cdots + \alpha_{k}\right)\right\}\\
&\cup \left\{\left(\alpha_{0} + 2\alpha_{1} + \cdots + 2\alpha_{j} + \alpha_{j+1} + \cdots + \alpha_{i} , \alpha_{i+1} + \cdots + \alpha_{k}\right) | k>i>j\right\}\\
&\cup \left\{\left(\alpha_{0} + 2\alpha_{1} + \cdots + 2\alpha_{i} + \alpha_{i+1} + \cdots + \alpha_{j} , \alpha_{i+1} + \cdots + \alpha_{k}\right) : k>j>i\right\}\\
&\cup\left\{\left(\alpha_{i+1} + \cdots + \alpha_{j} , \alpha_{0} + 2\alpha_{1} + \cdots + 2\alpha_{i} + \alpha_{i+1} + \cdots + \alpha_{k}\right) : k>j>i\right\}
\end{align*}
Therefore,
\begin{align*}
l\left(\beta\right) =& \min\left\{[0 , 1 , \ldots , j , 1 , \ldots , k]\right\} \\
&\cup \left\{[j , \ldots , 1 , 0 , 1 , \ldots , k]\right\} \\
&\cup \left\{[i , \ldots , 1 , 0 , 1 , \ldots , j , i+1 , \ldots , k] | j>i>k\right\} \\
&\cup \left\{[i+1 , \ldots , j , i , \ldots , 1 , 0 , 1 , \ldots , k]| j>i>k\right\} \\
=& [j , \ldots , 1 , 0 , 1 , \ldots , k].
\end{align*}

\item Arguing as in the type $A$ case gives $l(\af_i+\cdots+\af_j)=[i,\ldots,j]$ for $1\leq i\leq j$.

Observe that the remaining roots may be written as $\bt=\af_0+\cdots+\af_j+\af_2+\cdots+\af_k$ for $0\leq j<k$ and $k\geq2$. For the base case we have that $l(\af_0+\af_2)=[0,2]$.

Now, let $\bt=\af_0+\cdots+\af_j+\af_2+\cdots+\af_k$, $0\leq j<k$, $k\geq2$ (assume $k>2$ if $j=0$). We may assume by induction that if either $i<j$ and $m\leq k$, or $i=j$ and $m<k$ that
$$l(\af_0+\cdots+\af_i+\af_2+\cdots+\af_m)=[i,\ldots 1,0,2,\ldots,m].$$
We have
\begin{align*}
C\left(\beta \right) =& \left\{\left(\alpha_{0} + \alpha_{1} + \cdots + \alpha_{j} + \alpha_{2} + \cdots + \alpha_{i} , \alpha_{i+1} + \cdots + \alpha_{k}\right) : k>i>j , i \geq 2\right\} \\
&\cup \left\{\left(\alpha_{0} + \alpha_{1} + \cdots + \alpha_{i} + \alpha_{2} + \cdots + \alpha_{j} , \alpha_{i+1} + \cdots + \alpha_{k}\right) : k>j>i\geq0 , j \geq 2\right\}  \\
&\cup \left\{\left(\alpha_{i+1} + \cdots + \alpha_{j} , \alpha_{0} + \alpha_{1} + \cdots + \alpha_{i} + \alpha_{2} + \cdots + \alpha_{k}\right) : k>j>i\geq0\right\}
\end{align*}
Thus,
\begin{align*}
l\left(\beta\right) =& \min\left\{[j , \ldots , 1 , 0 , 2 , \ldots , k]\right\} \\
&\cup \left\{[i , \ldots , 1 , 0 , 2 , \ldots , j , i+1 , \ldots , k] | k>j>i\geq0 , j \geq 2\right\} \\
&\cup \left\{[i+1 , \ldots , j , i , \ldots , 1 , 0 , 2 , \ldots , k] | k>j>i\geq0\right\} \\
=& [j , \ldots , 1 , 0 , 2 , \ldots , k]
\end{align*}
\end{enumerate}
\end{pff}

\subsection{Root Vectors}\label{A:Root Vectors}

\begin{prp} In type $A_r$,
\[
b^*[i,\ldots,j]=[i,\ldots,j],\;\;\;0\leq i\leq j<r.
\]
\end{prp}

\begin{pff}
Proceed by induction on $j-i$, the case $j-i=0$ being trivial. Assume that
$i<j$ and $r_{[i+1,\ldots,j]}=(q-q^{-1})^{j-i-1}[i+1,\ldots,j]$. Using equation
\eqref{E:qShuffle} we deduce that
\begin{align*}
r_{[i,\ldots,j]}=\Xi(\la[i,\ldots,j]\ra)&=\Xi([[i],\la[i+1,\ldots,j]\ra]_q)\\
    &=\Xi([i])*\Xi(\la[i+1,\ldots,j]\ra)-q^{-1}\Xi(\la[i+1,\ldots,j]\ra)*\Xi([i])\\
    &=[i]*r_{[i+1,\ldots,j]}-q^{-1}r_{[i+1,\ldots,j]}*[i]\\
    &=(q-q^{-1})^{j-i-1}[i]*[i+1,\ldots,j]-[i]\overline{*}[i+1,\ldots,j]\\
    &=(q-q^{-1})^{j-i-1}([i]*[i+1]-[i]\overline{*}[i+1])[i+2,\ldots,j]\\
    &=(q-q^{-1})^{j-i-1}(q[i,i+1]-q^{-1}[i,i+1])[i+2,\ldots,j]\\
    &=(q-q^{-1})^{j-i}[i,\ldots,j].
\end{align*}
Finally, using \ref{E:Proportionality} we deduce that $b^*[i,\ldots,j]=E^*[i,\ldots,j]=[i,\ldots,j]$.
\end{pff}

\begin{prp} In type $B_r$:
\begin{align*}
b^*[i,\ldots,j]&=[i,\ldots,j],\;\;\;0\leq i\leq j<r\\
b^*[j,\ldots,0,0,\ldots,k]&=[2]_0[j,\ldots,0,0,\ldots,k],\;\;\;0\leq j<k<r.
\end{align*}
\end{prp}

\begin{pff}
The first formula follows easily by induction on $j-i$ as in the type $A$ case. We prove the second formula by induction on $j$ and $k$ with $j<k$, using \eqref{E:inductiveqshuffle}, \eqref{E:qbracket}, and \eqref{E:Lyndonbracketing} for the computations.

Observe that for $k\geq1$, $r_{[0,1,\ldots,k]}=(q^2-q^{-2})^k[0,1,\ldots,k]$,
which can be proved easily by downward induction on $j$, $0\leq j<k$, using
\eqref{E:inductiveqshuffle} and
\[
r_{[j,\ldots,k]}=\Xi(\la[j,\ldots,k]\ra)=\Xi([[j],\la[j+1,\ldots,k]\ra]_q)=[j]*r_{[j+1,\ldots,k]}-q^{-2}r_{[j+1,\ldots,k]}*[j].
\]
By \eqref{E:inductiveqshuffle}, we have
\begin{align*}
[0]*[0,1]-[0,1]*[0]&=[0,1,0]+q^2([0]*[0])\cdot[1]-([0]*[0])\cdot[1]-[0,1,0]\\
&=(q^2-1)([0,0]+q^{-2}[0,0])\cdot[1]=(q^2-q^{-2})\cdot[0,0,1]
\end{align*}
Therefore, applying \eqref{E:Lyndonbracketing} and the relevant definitions, we deduce that
\begin{align*}
r_{[0,0,1]}&=\Xi(\la[0,0,1]\ra)\\
    &=\Xi([[0],\la[0,1]\ra]_q)\\
    &=[0]*r_{[0,1]}-r_{[0,1]}*[0]\\
    &=(q^2-q^{-2})([0]*[0,1]-[0,1]*[0])\\
    &=(q^2-q^{-2})^2[0,0,1]
\end{align*}
Once again, using \eqref{E:inductiveqshuffle}, we deduce that for all $k\geq2$,
\begin{align}\label{E:DblSegReduction0}
[0]*[0,\ldots,k]-[0,\ldots,k]*[0]=([0]*[0,\ldots,k-1]-[0,\ldots,k-1]*[0])\cdot[k].
\end{align}
Assume $k\geq2$. Then, $(\af_0,\af_0+\cdots+\af_k)=0$, so iterated applications
of \eqref{E:DblSegReduction0} yields
\begin{align*}
r_{[0,0,\ldots,k]}&=[0]*r_{[0,\ldots,k]}-r_{[0,\ldots,k]}*[0]\\ &=(q^2-q^{-2})^k([0]*[0,\ldots,k]-[0,\ldots,k]*[0])\\
    &=(q^2-q^{-2})^k([0]*[0,1]-[0,1]*[0])\cdot[2,\ldots, k]\\
    &=(q^2-q^{-2})^{k+1}[0,0,\ldots, k]
\end{align*}

Now, assume that $k\geq 2$, and $0<j<k$. To compute $r_{[j,\ldots,0,0,\ldots,k]}$, we need the following. For $|j-k|>1$,
\begin{align}\label{E:DblSegReduction1}
[j]*[j-1,\ldots,&k]-q^{-2}[j-1,\ldots,k]*[j]\\
    \nonumber&=([j]*[j-1,\ldots,k-1]-q^{-2}[j-1,\ldots,k-1]*[j])\cdot[k].
\end{align}
For $j=k-1$,
\begin{align}\label{E:DblSegReduction2}
[j]*[j-1,\ldots,0,&0,\ldots,j+1]-q^{-2}[j-1,\ldots,0,0,\ldots,j+1]*[j]\\
    \nonumber&=(q^2[j]*[j-1,\ldots,0,0,\ldots,j]-q^{-2}[j-1,\ldots,0,0,\ldots,j]*[j])\cdot[j+1].
\end{align}
Finally,
\begin{align}\label{E:DblSegReduction3}
q^2[j]*[j-1,\ldots,0,&0,\ldots,j]-q^{-2}[j-1,\ldots,0,0,\ldots,j]*[j]\\
    \nonumber&=([j]*[j-1,\ldots,0,0,\ldots,j-2]-q^{-2}[j-1,\ldots,0,0,\ldots,j-2]*[j])\cdot[j,j+1].
\end{align}
Indeed, \eqref{E:DblSegReduction1} and \eqref{E:DblSegReduction2} are
straightforward applications of \eqref{E:inductiveqshuffle}. Equation
\eqref{E:DblSegReduction3} involves a little more calculation:
\begin{align*}
q^2[j]*[j-1,\ldots,0,&0,\ldots,j]-q^{-2}[j-1,\ldots,0,0,\ldots,j]*[j]\\
=&q^2[j-1,\ldots,0,0,\ldots,j,j]+q^{-2}([j]*[j-1,\ldots,0,0,\ldots,j-1]\\
    &-[j-1,\ldots,0,0,\ldots,j-1]*[j])\cdot[j]-q^{-2}[j-1,\ldots,0,0,\ldots,j,j]\\
    =&(q^2-q^{-2})\cdot[j-1,\ldots,0,0,\ldots,j,j]+q^{-2}([j-1,\ldots,0,0,\ldots,j]\\&+q^2([j]*[j-1,\ldots,0,0,\ldots,j-2])\cdot[j-1]
    -([j-1,\ldots,0,0,\ldots,j-2]*[j])\cdot[j-1]\\&-q^4[j-1,\ldots,0,0,\ldots,j])\cdot[j]\\
    =&([j]*[j-1,\ldots,0,0,\ldots,j-2]-q^{-2}[j-1,\ldots,0,0,\ldots,j-2]*[j])\cdot[j,j+1],
\end{align*}
Note that \eqref{E:DblSegReduction1} holds for both $[j-1,j,\ldots,k]$ and $[j-1,\ldots,0,0,\ldots,k]$.

Now, assume that we have shown that
$r_{[j-1,\ldots,0,0,\ldots,k]}=(q^2-q^{-2})^{j+k}[j-1,\ldots,0,0,\ldots,k]$.
Then, since $(\af_j,2\af_0+\cdots+2\af_{j-1}+\af_j+\cdots+\af_k)=-2$,
\begin{align*}
r_{[j,\ldots,0,0,\ldots,k]}=&[j]*r_{[j-1,\ldots,0,0,\ldots,k]}-q^{-2}r_{[j-1,\ldots,0,0,\ldots,k]}*[j]&\\
    =&(q^2-q^{-2})^{j+k}[j]*[j-1,\ldots,0,0,\ldots,k]-q^{-2}[j-1,\ldots,0,0,\ldots,k]*[j]&\\
    =&(q^2-q^{-2})^{j+k}([j]*[j-1,\ldots,0,0,\ldots,j+1]\\&-q^{-2}[j-1,\ldots,0,0,\ldots,j+1]*[j])\cdot[j+2,\ldots,k]
     &\mbox{by \eqref{E:DblSegReduction1}}\\
    =&(q^2-q^{-2})^{j+k}(q^2[j]*[j-1,\ldots,0,0,\ldots,j]&\\&-q^{-2}[j-1,\ldots,0,0,\ldots,j]*[j])\cdot[j+1,\ldots,k]
     &\mbox{by \eqref{E:DblSegReduction2}}\\
    =&(q^2-q^{-2})^{j+k}([j]*[j-1,\ldots,0,0,\ldots,j-2]&\\&-q^{-2}[j-1,\ldots,0,0,\ldots,j-2]*[j])\cdot[j,\ldots,k]
     &\mbox{by \eqref{E:DblSegReduction3}}\\
    =&(q^2-q^{-2})^{j+k}([j]*[j-1]-q^{-2}[j-1]*[j])\cdot[j-2,\ldots,0,0,\ldots,k]&\mbox{by \eqref{E:DblSegReduction1}}\\
    =&(q^2-q^{-2})^{j+k+1}[j,\ldots,0,0,\ldots,k].&
\end{align*}

Finally, one computes using \eqref{E:KashiwaraForm2} and \eqref{E:Nfcn} that
the coefficient of $[j,\ldots,0,0,\ldots,k]$ in \eqref{E:Proportionality} is
\[
\frac{(-1)^{j+k}(q^2-q^{-2})^{j+k+1}(1-q^4)}{q^{-2(j+k)}(1-q^2)^2(1-q^4)^{j+k}}=[2]_0^2,
\]
so the result follows.
\end{pff}

\begin{prp}\label{P:C root vectors} In type $C_r$:
\begin{align*}
b^*[i,\ldots,j]&=[i,\ldots,j],\;\;\;0\leq i\leq j<r,\\
b^*[j,\ldots,1,0,1,\ldots,k]&=[j,\ldots,1,0,1,\ldots,k],\;\;\;1\leq j<k<r,\\
b^*[0,\ldots,j,1,\ldots,j]&=q[0]\cdot([1,\ldots,j]*[1,\ldots,j]),\;\;\;1\leq j<r.
\end{align*}
\end{prp}

\begin{pff}
The first to formulae can be proved by induction as in the type $A$ case.

We now prove that $b^*[0,\ldots,j,1\ldots,j]=q[0]\cdot([1,\ldots,j]*[1,\ldots,j])$. Our argument is essentially the same as \cite[Lemma 53]{lec}. Indeed, $[1,\ldots,j]$ belongs to $\U_q$, so $[1,\ldots,j]*[1,\ldots,j]$ belongs to $\U_q$. Using \cite[Theorem 5]{lec}, we deduce that $f=[0]\cdot([1,\ldots,j]*[1,\ldots,j])$ belongs to $\U_q$. Clearly $\min(f)=[0,\ldots,j,1,\ldots,j]$ so by Theorem \ref{T:minbg}(i) $f$ is proportional to $b^*[0,\ldots,j,1,\ldots,j]$. Finally, using Theorem \ref{T:minbg}(i) and \eqref{E:Normalized Eg} with $g=[1,\ldots,j,1,\ldots,j]$ we obtain the result.
\end{pff}

\begin{prp} In type $D_r$:
\begin{align*}
b^*[0]&=[0]\\
b^*[0,2,\ldots,i]&=[0,2,\ldots,i],\;\;\;2\leq i<r,\\
b^*[i,\ldots,j]&=[i,\ldots,j],\;\;\;1\leq i\leq j<r,\\
b^*_{[0,\ldots,j]}&=[1,0,2,\ldots,j]+[0,1,2,\ldots,j],\;\;\;2\leq j<r,\\
b^*[j,\ldots,2,1,0,2,\ldots,k]&=[j,\ldots,2,1,0,2,\ldots,k]+[j,\ldots,2,0,1,2,\ldots,k],\;\;\;2\leq
j<k<r.
\end{align*}
\end{prp}

\begin{pff}
All cases follow by an easy induction argument that we leave as an exercise for the reader.
\end{pff}

\subsection{The $F_4$ case}\label{A:F} Calculations available upon request.

\subsection{The $E_8$ case}\label{A:E8}

\begin{center}
\begin{tabular}{|c|c|}\hline
Height&Good Lyndon Words\\
\hline\hline
1&$[0], [1], [2],[3], [4], [5], [6], [7]$\\
2&$[01],[13],[23], [34], [45], [56], [67]$\\
3&$[013], [123], [134], [234], [345], [456], [567]$\\
4&$[0123],[2134],[2345],[0234],[1345],[3456],[4567]$\\
5&$[02134],[21345],[13456],[32134],[02345],[23456],[34567]$\\
6&$[021345],[213456],[302134],[321345],[023456],[134567],[234567]$\\
7&$[2134567],[2302134],[4321345],0213456],[0234567],[3213456]$\\
8&$[02134567],[23021345],[43012345],[32134567],[30213456],[43213456]$\\
9&$[423021345],[230213456],[43213456],[543213456],[430213456],[423021345],[302134567]$\\
10&$[3423021345],[4230213456],[5430213456],[5432134567],[2302134567],[4302134567]$\\
11&$[13423021345],[54230213456],[34230213456],[42302134567],[54302134567],[5432134567]$\\
12&$[534230213456],[134230213456],342302134567],[54230213456],[654302134567]$\\
13&$[4534230213456],[5134230213456],[1342302134567],[5342302134567],[6542302134567]$\\
14&[45134230213456],[51342302134567],[45342302134567],[65342302134567]\\
15&[314534230213456],[451342302134567],[651342302134567],[645342302134567]\\
16&[2314534230213456],[3145342302134567],[6451342302134567],[5645342302134567]\\
17&[02314534230213456],[23145342302134567],[63145342302134567],[56451342302134567]\\
18&[023145342302134567],[623145342302134567],[563145342302134567]\\
19&[6023145342302134567],[562345342302134567],[4563145342302134567]\\
20&[56023145342302134567],[4562345342302134567]\\
21&[34562345342302134567],[456023145342302134567]\\
22&[134562345342302134567],[3456023145342302134567]\\
23&[13456023145342302134567],[23456023145342302134567]\\
24&[213456023145342302134567]\\
25&[3213456023145342302134567]\\
26&[43213456023145342302134567]\\
27&[543213456023145342302134567]\\
28&[6543213456023145342302134567]\\
29&[53423021345676451342302134567]\\\hline
\end{tabular}
\end{center}

\end{document}